\documentclass[12pt,leqno]{article}

=155mm\textheight =220mm

\usepackage{amsmath}
\usepackage{amscd}
\usepackage{amsmath}
\usepackage{amssymb}
\usepackage{latexsym}
\makeatletter

\@addtoreset{equation}{section}
\makeatother

\newtheorem{thm}{Theorem}[subsection]
\newtheorem{pr}[thm]{Proposition}
\newtheorem{df}[thm]{Definition}
\newtheorem{lm}[thm]{Lemma}
\newtheorem{cor}[thm]{Corollary}

\newtheorem{rmk}[thm]{Remark}
\newtheorem{eg}[thm]{Example}

\newcommand{\qed}{\hfill{\rule{7pt}{7pt}}}

\input xy
\xyoption{all}
\CompileMatrices

\begin{document}

\title{The characteristic class and ramification
of an $\ell$-adic \'etale sheaf}
\author{Ahmed Abbes and Takeshi Saito}
\maketitle

\hfill{\it Dedicated to Luc Illusie,
with admiration}

\begin{abstract}
We introduce the characteristic class of an $\ell$-adic \'etale sheaf
using a cohomological pairing due to Verdier (SGA5). 
As a consequence of the Lefschetz-Verdier 
trace formula, its trace 
computes the Euler-Poincar\'e characteristic of the sheaf.
We compare the characteristic class to
two other invariants arising from ramification theory.
One is the Swan class of Kato-Saito \cite{KS}
and the other is the 0-cycle class
defined by Kato 
for rank 1 sheaves in \cite{Kato2}.
\end{abstract}

\setcounter{section}0

Let $k$ be a perfect field, 
$\ell$ be a prime number invertible in $k$
and $\Lambda$ be a finite extension of either ${\mathbb F}_\ell$ or 
${\mathbb Q}_\ell$. Let $X$ be a separated $k$-scheme of finite type.
Let $f\colon X\to {\rm Spec}\ k$ denote the structural morphism 
and put ${\cal K}_X=Rf^!\Lambda$. 
For an 
\'etale sheaf $\cal F$
of $\Lambda$-modules on $X$,
its characteristic class
is defined as follows.
We put ${\cal H}=
R{\cal H}om({\rm pr}_2^*{\cal F},R{\rm pr}_1^!{\cal F})$
and ${\cal H}^*=
R{\cal H}om({\rm pr}_1^*{\cal F},R{\rm pr}_2^!{\cal F})$
on $X\times X$
and let $\Delta=X\subset X\times X$
denote the diagonal.
Then, we have
$1\in {\rm End}({\cal F})=
H^0_\Delta(X\times X,{\cal H})=
H^0_\Delta(X\times X,{\cal H}^*)$.
The natural pairing
${\cal H}\otimes{\cal H}^*
\to {\cal K}_{X\times X}$ induces
a pairing
$\langle\ ,\ \rangle\colon 
H^0_\Delta(X\times X,{\cal H})\otimes
H^0_\Delta(X\times X,{\cal H}^*)\to
H^0_\Delta(X\times X,{\cal K}_{X\times X})
=
H^0(X,{\cal K}_X)$.
We define the characteristic class of $\cal F$, denoted 
$C({\cal F})$, to be the pairing
$\langle1 ,1 \rangle
\in H^0(X,{\cal K}_X)$.
If $X$ is smooth of dimension $d$,
we have $C({\cal F})\in H^{2d}(X,\Lambda(d))$.
By the Lefschetz-Verdier trace formula \cite{SGA5}
Th\'eor\`eme 4.4,
the trace ${\rm Tr}\ C({\cal F})$
gives the Euler-Poincar\'e characteristic
$\chi(X_{\bar k},{\cal F})$
if $X$ is proper,
where $\bar k$ denotes a 
separable closure of $k$.

By devissage,
computations of the characteristic classes
are reduced to a computation of
$C(j_!{\cal F})$ 
where $j\colon U\to X$
is an open immersion,
$U$ is smooth and
${\cal F}$ is a smooth sheaf on $U$.
In this paper,
we compute the characteristic class of $C(j_!{\cal F})$
or rather the difference
$C(j_!{\cal F})-{\rm rank}\ {\cal F}\cdot C(j_!\Lambda)$
in terms of the ramification of ${\cal F}$
along the boundary $X\setminus U$.
More precisely, we prove that 
$C(j_!{\cal F})-{\rm rank}\ {\cal F}\cdot C(j_!\Lambda)$
is equal to the following two invariants,
under certain assumptions.
One is the Swan class ${\rm Sw}({\cal F})$
of $\cal F$ defined in \cite{KS}.
The other is 
the 0-cycle class $c_{\cal F}$ defined by Kato \cite{Kato2}
in rank 1 case.
We also define a localization
of the difference
$C(j_!{\cal F})-{\rm rank}\ {\cal F}\cdot C(j_!\Lambda)$
in $H^0_{X\setminus U}(X,{\cal K}_X)$
as a cohomology class with support
on the boundary.

The relation with
the Swan class
is a refinement of the generalized 
Grothendieck-Ogg-Shafarevich formula proved in 
\cite{KS}.
The proof uses a finite \'etale
covering of $U$ trivializing
the reduction modulo $\ell$ of 
the $\ell$-adic sheaf ${\cal F}$.
It is similar to that of
the generalized 
Grothendieck-Ogg-Shafarevich formula in loc.\ cit.
The key ingredients in the proof
are the compatibility of
the characteristic class with pull-back
and an explicit computation
of the characteristic class
in the tamely ramified case.
This argument works in arbitrary dimension
and in arbitrary rank,
under a certain mild assumption that
the sheaf is ``potentially of Kummer type''
(see Definition \ref{dfKum}).
However, since we use the pull-back to a covering,
we can not avoid a denominator.

The relation
with the other invariant
$c_{\cal F}$ in rank 1 case
is proved using a blow-up
at the ramification locus in the diagonal
$X\subset X\times X$.
This argument does not involve a denominator
and we get an integral result.
However, even to formulate the statement
at least for the moment,
we need to assume that the sheaf
has rank one.
Using this approach,
we obtain a new proof
of the Grothendieck-Ogg-Shafarevich
formula for curves 
without using
a covering trivializing the sheaf
or the Weil formula.
A crucial fact, used in the proof, is 
the following geometric interpretation
of the ramification theory
for Artin-Schreier-Witt characters
of Kato \cite{Kato},
proved in \cite{aml}. 

Let $(X\times X)'\to
X\times X$ be the log blow-up
and $X\to (X\times X)'$
be the log diagonal map.
We consider the blow-up
$(X\times X)''\to (X\times X)'$
at the Swan divisor
$D_{\cal F}\subset X$
regarded as a closed subscheme
of $(X\times X)'$
by the log diagonal map.
Then, the sheaf 
${\cal H}om({\rm pr}_2^*{\cal F},{\rm pr}_1^*{\cal F})$
on $U\times U$
is unramified along the
exceptional divisor of
$(X\times X)''$.
Further the restriction
of a smooth extension
of ${\cal H}om({\rm pr}_2^*{\cal F},{\rm pr}_1^*{\cal F})$
on the exceptional divisor
is the Artin-Schreier sheaf
defined by the refined Swan conductor.
This interpretation means that
one can kill the ramification of a sheaf
not only by a ramified covering
but also by blowing-up the
diagonal on the ramification locus.
More quantitatively,
the Swan conductor measures 
the necessary blow-up
to kill the ramification.
This fits nicely with the ramification theory
in \cite{AS}.
The authors expect that this argument
should work with arbitrary rank
(cf. \cite{aml}).

The paper is organized as follows.
In Section 1, we recall some sorites on cohomological correspondences
and the Verdier pairing, and we give some complements.
In Section 2, we define the characteristic class of a cohomological 
correspondence. Proposition (\ref{prpink}) gives a recipe to compute it.  
It plays a crucial role in this paper.
Then we compute the characterstic class under some tameness condition.
In Section 3, we introduce the notion of a {\em sheaf of Kummer type}, 
and recall the definition of the Swan class of Kato-Saito \cite{KS}. 
For a sheaf potentially of Kummer type, we prove the main formula 
of this paper relating its characteristic class and its Swan class. 
In Section 4, we study sheaves of rank one.
We recall the ramification theory
of Artin-Schreier-Witt characters in \cite{Kato}. 
Then, we prove an integral relation between the characteristic class and the
0-cycle class $c_{\cal F}$ for a rank one sheaf, by blowing-up the diagonal.
In Section 5,
we define a 
cohomology class 
$C_{X\setminus U}^0({\cal F})
\in H^0_{X\setminus U}(X,{\cal K}_X)$
with support on the boundary
and show that it is a refinement 
of the difference
$C(j_!{\cal F})-
{\rm rank}\ {\cal F}\cdot C(j_!\Lambda)$.

The authors thank Luc Illusie and Kazuya Kato
for inspiring discussions.
They thank Illusie for pointing out an error
in Lemma \ref{lmadd} in an earlier version
and informing the reference \cite{fer}.

\bigskip
\noindent
{\bf \Large Notation}
\smallskip

In this paper, $k$ denotes a field.
In sections 3 and 4,
we will assume $k$ perfect.
Schemes over $k$
are assumed separated and
of finite type.
Thus the diagonal map
$\delta\colon X=\Delta_X\to X\times X$
is a closed immersion.
For a divisor with simple normal crossings
of a smooth scheme over $k$,
we assume that
the irreducible components 
and their intersections
are also smooth over $k$.

The letter $\ell$ denotes
a prime number invertible in $k$
and $\Lambda$ denotes
a finite commutative ${\mathbb Z}_\ell$-algebra.
For a scheme $X$ over $k$,
$D_{\rm ctf}(X)$
denotes the derived category of
complexes of $\Lambda$-modules
of finite tor-dimension on
the \'etale site of $X$
with constructible cohomology
\cite{Rapport} 4.6.
We omit to write 
$R$ or $L$
to denote the derived functors
unless otherwise
stated explicitly
or for $R{\cal H}om$.
Let ${\cal K}_X$ denote $f^!\Lambda$
where $f\colon X\to {\rm Spec}\ k$
is the structure map and
let ${\mathbf D}$ denote the functor
$R{\cal H}om(\ ,{\cal K}_X)$.
For objects ${\cal F}$
and ${\cal G}$
of $D_{\rm ctf}(X)$
and  $D_{\rm ctf}(Y)$
on schemes $X$ and $Y$ over $k$,
${\cal F}\boxtimes{\cal G}$
denotes
${\rm pr}_1^*{\cal F}\otimes
{\rm pr}_2^*{\cal G}$
on $X\times Y$.

\section{Preliminaries on cohomological correspondences.}

We recall some generalities on \'etale cohomology
and constructions in \cite{SGA5}.

\subsection{Review on base change maps,
cycle class maps etc.}

Let 
$$\begin{CD}
X'@>{f'}>> Y'\\
@V{g'}VV @VVgV\\
X@>f>> Y
\end{CD}$$
be a commutative diagram 
of schemes over $k$.
We have base change
morphisms
of functors
\begin{eqnarray}
g^*f_*&\to& f'_*g^{\prime *},
\label{g^*f_*}\\
f'_!g^{\prime !}&\to& g^!f_!
\label{f^!g_!}
\end{eqnarray}
(\cite{4XVII} 4.1.5,
\cite{4XVIII} 3.1.13.2).
If $f$ is proper, 
if $g$ is an open immersion
and if the square is cartesian,
the maps 
(\ref{g^*f_*}) and (\ref{f^!g_!})
are isomorphisms and 
are inverse to each other.
If $g'$ is proper,
the base change map
\begin{equation}
g^*f_!\to f'_!g^{\prime *}
\label{g^*f_!}
\end{equation}
is similarly defined as the adjoint of
the composition
$f_!
\to f_!g'_*g^{\prime *}
= f_!g'_!g^{\prime *}
\to g_!f'_!g^{\prime *}
\to g_*f'_!g^{\prime *}$.
The base change maps (\ref{g^*f_*})
and (\ref{g^*f_!}) form
a commutative diagram
$$\begin{CD}
g^*f_!@>>> f'_!g^{\prime *}\\
@VVV @VVV\\
g^*f_*@>>> f'_*g^{\prime *}.
\end{CD}$$

For a commutative diagram
$$\begin{CD}
X'@>{f'}>> Y'@>{g'}>> Z'\\
@V{h''}VV @V{h'}VV @VV{h}V\\
X@>f>> Y@>g>> Z
\end{CD}$$
of schemes over $k$,
the base change morphisms (\ref{f^!g_!})
are transitive.
Namely, we have a commutative diagram
\begin{equation}
\xymatrix{
g'_!f'_!h^{\prime\prime!}\ar[rr]\ar[dr]&&
g'_!h^{\prime !}f_!\ar[rr]&&
h^!g_!f_!\ar[dl]\\
&(g'f')_!h^{\prime\prime!}\ar[rr]&&
h^!(gf)_!}
\label{f^!g_!2}
\end{equation}
where the slant arrows
are induced by
the composition isomorphisms.
We also have a similar commutative diagram
for base change morphisms (\ref{g^*f_*}).
If $h'$ and $h''$ are proper,
we have a similar commutative diagram
for (\ref{g^*f_!}).

For a morphism $f\colon X\to Y$
of schemes over $k$
and objects ${\cal F}$
and ${\cal G}$
of $D_{\rm ctf}(Y)$,
we recall the definition
of the canonical isomorphism
\begin{equation}
R{\cal H}om(f^*{\cal F},f^!{\cal G})
\to f^!R{\cal H}om({\cal F} ,{\cal G})
\label{fHom}
\end{equation} 
(\cite{4XVIII} 3.1.12.2).
A canonical map
${\cal F}\otimes
f_!R{\cal H}om(f^*{\cal F},f^!{\cal G})
\to 
f_!(f^*{\cal F}\otimes
R{\cal H}om(f^*{\cal F},f^!{\cal G}))$
is defined in
(\cite{4XVIII} 3.1.11.2).
By composing
$f_!(f^*{\cal F}\otimes
R{\cal H}om(f^*{\cal F},f^!{\cal G}))
\to f_!f^!{\cal G}\to {\cal G}$,
we obtain
${\cal F}\otimes
f_!R{\cal H}om(f^*{\cal F},f^!{\cal G})
\to {\cal G}$.
This induces the map (\ref{fHom})
by adjunction.
If ${\cal F}$ and ${\cal G}$
are objects of the filtered derived
category \cite{cxct} Chap.\ V \S\S1-3,
the isomorphism (\ref{fHom})
is an isomorphism in
the filtered derived category.
Recall that
an object of $DF_{\rm ctf}(X)$
is a complex ${\cal F}$
of $\Lambda$-modules
on the \'etale site of $X$
with finite filtration $F_\bullet$
such that $Gr^F_\bullet {\cal F}$
is an object of
$D_{\rm ctf}(X)$.
In particular,
we have an isomorphism
\begin{equation}
Hom_{DF_{\rm ctf}(X)}(f^*{\cal F},f^!{\cal G})
\to \Gamma(X,
f^!F_0R{\cal H}om({\cal F} ,{\cal G})).
\label{fHomg}
\end{equation}

For a morphism $f\colon X\to Y$
of schemes over $k$
and objects ${\cal F}$
and ${\cal G}$
of $D_{\rm ctf}(Y)$,
a canonical map
\begin{equation}
f^*{\cal F} \otimes f^!{\cal G}
\to f^!({\cal F} \otimes {\cal G})
\label{FoxG}
\end{equation} 
is defined as follows.
We have the canonical isomorphism
$f_!(f^*{\cal F} \otimes f^!{\cal G})
\to 
{\cal F} \otimes f_!f^!{\cal G}$
of the projection formula
(\cite{4XVII} 5.2.9).
By composing with the adjunction map
$f_!f^!{\cal G}\to {\cal G}$,
we obtain a map
$f_!(f^*{\cal F} \otimes f^!{\cal G})
\to 
{\cal F} \otimes {\cal G}$.
Thus, the canonical map
(\ref{FoxG})
is defined by taking the adjoint.

For a flat morphism
$f\colon X\to Y$ of fiber dimension $d$
of schemes over $k$,
the trace map
$f_!\Lambda(d)[2d]\to \Lambda$
is defined (\cite{4XVIII} 2.9)
and the canonical class map
\begin{equation}
\Lambda(d)[2d]\to Rf^!\Lambda
\label{eqcan}
\end{equation}
is defined as its adjoint
(\cite{4XVIII} 3.2.3).
If $f$ is smooth,
the map (\ref{eqcan})
is an isomorphism
(\cite{4XVIII} 3.2.5).

Let
$f\colon X\to Y$ 
be a morphism of
smooth schemes
of the same dimension $n$ over $k$.
The map $f$ is the composition
${\rm pr_2}\circ \gamma$
of the graph map $\gamma\colon 
X\to X\times Y$
and the second projection ${\rm pr}_2\colon 
X\times Y\to Y$.
Since $\gamma$ is a section 
of the first projection
${\rm pr}_1\colon 
X\times Y\to X$
and the projections
are smooth of relative dimension $n$,
a canonical isomorphism
$\Lambda\to
f^!\Lambda$
is defined as the composition 
$\Lambda
\to
\gamma^!{\rm pr}_1^!\Lambda
\overset{\sim} \gets 
\gamma^!\Lambda(n)[2n]
\to 
\gamma^!{\rm pr}_2^!\Lambda
\to 
f^!\Lambda$.
Thus the canonical map 
(\ref{FoxG}) for ${\cal G}=\Lambda$
induces a functor
\begin{equation}
f^*\to f^!.
\label{f^*f^!}
\end{equation}

We recall the cycle class map.
Let $C$ be a scheme over $k$
and $d\ge 0$ be an integer.
For an integral closed subscheme $V\subset C$
of dimension $d$,
the canonical class map
$\Lambda(d)[2d]\to {\cal K}_V$ gives
the cycle class
$[V]\in H^0(V,{\cal K}_V(-d)[-2d])=H^{-2d}_V(C,{\cal K}_C(-d))$.
The cycle classes define a map
${\rm cl}\colon Z_d(C)\to H^{-2d}(C,{\cal K}_C(-d))$
where $Z_d(C)$
denotes the free abelian group 
generated by the 
integral closed subschemes of $C$
of dimension $d$.
The cycle class map
$Z_d(C)\to H^{-2d}(C,{\cal K}_C(-d))$
factors through the quotient
by rational equivalence
and induces a map
\begin{equation}
\begin{CD}
{\rm cl}_C\colon CH_d(C)@>>> H^{-2d}(C,{\cal K}_C(-d))
\end{CD}
\label{eqcycleCH}
\end{equation}

If $C$ is a closed subscheme of
a smooth scheme $X$ of dimension $d+c$,
the target $H^{-2d}(C,{\cal K}_C(-d))$
is identified with
$H^{-2d}_C(X,{\cal K}_X(-d))=
H^{2c}_C(X,\Lambda(c)).$
The cycle class map
and the intersection product
are compatible in the following sense.
Let $X$ and $Y$ be smooth schemes
of dimension $d+c$ and $e+c$ over $k$
and $f\colon X\to Y$ be a morphism over $k$.
Let $C\subset X$ and $D\subset Y$ be closed
subschemes
satisfying $f^{-1}(D)\subset C$ set-theoretically.
Then the diagram
\begin{equation}
\begin{CD}
CH_e(D)@>{{\rm cl}_D}>> H^{2c}_D(Y,\Lambda(c))\\
@V{f^!}VV @VV{f^*}V\\
CH_d(C)@>{{\rm cl}_C}>> H^{2c}_C(X,\Lambda(c))
\end{CD}
\label{eqcycpull}
\end{equation}
is commutative by \cite{KS} Lemma 2.1.2.
The left vertical arrow
$f^!\colon CH_e(D)\to CH_d(C)$
denotes the Gysin map in the intersection theory
\cite{fulton}.
The commutative diagram (\ref{eqcycpull})
for the diagonal map
$\delta\colon X\to Y=X\times X$
implies the compatibility
\begin{equation}
{\rm cl}(V)\cup {\rm cl}(W)=
{\rm cl}(V,W)_X
\label{eqcupint}
\end{equation}
in $H^{2(c+d)}_{V\cap W}(X,\Lambda (c+d))$
of the cup-product and
the intersection theory
for closed subschemes $V$ and $W$ in $X$
of codimensions $c$ and $d$.

\subsection{Cohomological correspondences.}

We recall definitions on cohomological
correspondences
following \cite{SGA5}.
In order to remain compatible with the convention
in \cite{4XVII} 1.1.10 (ii),
we have chosen to switch the factors in \cite{SGA5}
and consider a cohomological
correspondence as a map
from the second factor to the first factor.

\begin{df}\label{dfcohcor}
Let $X$ and $Y$ be schemes
over $k$ and ${\cal F}$ and ${\cal G}$
be objects of
$D_{\rm ctf}(X)$
and of
$D_{\rm ctf}(Y)$
respectively.
We call a correspondence
between $X$ and $Y$
a scheme $C$ over $k$
and morphisms
$c_1\colon C\to X$ and
$c_2\colon C\to Y$ over $k$.
We put
$c=(c_1,c_2)\colon C\to X\times Y$
the corresponding morphism.
We call a morphism
$u\colon c_2^*{\cal G}\to c_1^!{\cal F}$
a cohomological correspondence from
${\cal G}$ to ${\cal F}$ on $C$.
\end{df}

The definition here is slightly different from
that given in \cite{SGA5} 3.2
if $c\colon C\to X\times Y$ is not proper.
If $c\colon C\to X\times Y$
is the graph map $\gamma\colon X\to X\times Y$
of a map $f\colon X\to Y$,
a cohomological correspondence from
${\cal G}$ to ${\cal F}$ on $C$
is nothing but a map
$f^*{\cal G}\to {\cal F}$.

By the isomorphism (\ref{fHom})
\begin{equation}
c^!R{\cal H}om(
{\rm pr}_2^*{\cal G},{\rm pr}_1^!{\cal F})
\to 
R{\cal H}om(
c_2^*{\cal G},c_1^!{\cal F})
\label{eqw!}
\end{equation}
and by the adjunction,
the following three notions
are equivalent:
\begin{enumerate}
\item[(1)]
a cohomological correspondence
$u\colon c_2^*{\cal G}\to c_1^!{\cal F}$
on $C$.
\item[(2)]
a section of 
$c^!R{\cal H}om(
{\rm pr}_2^*{\cal G},{\rm pr}_1^!{\cal F})$
on $C$,
or equivalently
a map
$\Lambda_C\to c^!R{\cal H}om(
{\rm pr}_2^*{\cal G},{\rm pr}_1^!{\cal F})$.
\item[(3)]
a map
$c_!\Lambda_C\to R{\cal H}om(
{\rm pr}_2^*{\cal G},{\rm pr}_1^!{\cal F})$.
\end{enumerate}
We will identify them freely in the following.
If ${\cal F}$ and ${\cal G}$
are objects in the filtered derived
categories,
a homomorphism
$u\colon c_2^*{\cal G}\to c_1^!{\cal F}$
in $DF_{\rm ctf}(C)$
corresponds bijectively to
a section of 
$c^!F_0R{\cal H}om(
{\rm pr}_2^*{\cal G},{\rm pr}_1^!{\cal F})$
by the isomorphism
(\ref{fHomg}).
If $c\colon C\to X\times Y$
is a closed immersion,
the three notions above
are further equivalent to the following:
\begin{enumerate}
\item[(4)]
an element of
$H^0_C(X\times Y,R{\cal H}om(
{\rm pr}_2^*{\cal G},{\rm pr}_1^!{\cal F}))$.
\end{enumerate}
Recall that a canonical
isomorphism
\begin{equation}
\begin{CD}
R{\cal H}om(
{\rm pr}_2^*{\cal G},{\rm pr}_1^!{\cal F})
@>>>{\cal F}
\boxtimes^L
{\mathbf D}{\cal G}\end{CD}
\label{eqFG}
\end{equation}
is defined in 
\cite{SGA5} (3.1.1).

A typical example of
a cohomological correspondence is
given as follows.
Assume $X$ and $Y$ are smooth of dimension $d$
over $k$
and $c=(c_1,c_2)\colon C\to X\times Y$
is a closed immersion.
Let ${\cal F}$ and ${\cal G}$
be sheaves of free $\Lambda$-modules on 
$X$ and $Y$ respectively
and assume ${\cal G}$ is smooth.
Then, the canonical map
$c^*{\cal H}om({\rm pr}_2^*{\cal G},
{\rm pr}_1^*{\cal F})
\to {\cal H}om(c_2^*{\cal G},c_1^*{\cal F})$
is an isomorphism
and we identify
$Hom(c_2^*{\cal G},c_1^*{\cal F})
=\Gamma(C,c^*{\cal H}om({\rm pr}_2^*{\cal G},
{\rm pr}_1^*{\cal F}))$.
Since ${\rm pr}_1\colon X\times Y\to X$ is smooth,
we have a canonical isomorphism
${\rm pr}_1^*{\cal F}(d)[2d]
\to {\rm pr}_1^!{\cal F}$
(\cite{4XVIII} 3.2.5)
and we identify
$R{\cal H}om({\rm pr}_2^*{\cal G},
{\rm pr}_1^!{\cal F})=
{\cal H}om({\rm pr}_2^*{\cal G},
{\rm pr}_1^*{\cal F})(d)[2d]$.
Then the cycle class map
$CH_d(C)\to H^{2d}_C(X\times Y,\Lambda(d))$
induces a pairing
\begin{equation}
\begin{CD}
CH_d(C)\otimes
Hom(c_2^*{\cal G},c_1^*{\cal F})
@>>> H^{2d}_C(X\times Y,\Lambda(d))
\otimes
\Gamma(C,c^*{\cal H}om({\rm pr}_2^*{\cal G},
{\rm pr}_1^*{\cal F}))\\
@>>> H^0_C(X\times Y,
R{\cal H}om({\rm pr}_2^*{\cal G},
{\rm pr}_1^!{\cal F})).
\end{CD}
\label{equgam}
\end{equation}
In other words,
the pair $(\Gamma,\gamma)$
of a cycle class
$\Gamma\in CH_d(C)$
and a homomorphism
$\gamma\colon c_2^*{\cal F}\to c_1^*{\cal F}$
defines a cohomological correspondence
$u(\Gamma,\gamma)$.

We recall the definition of
the push-forward
of a cohomological correspondence.
We consider a commutative diagram
\begin{equation}
\begin{CD}
X@<{c_1}<< C@>{c_2}>> Y\\
@VfVV @VhVV @VVgV\\
X'@<{c'_1}<< C'@>{c'_2}>> Y'
\end{CD}\label{eqpush}
\end{equation}
of schemes over $k$.
A canonical isomorphism
\begin{equation}
(f\times g)_*R{\cal H}om(
{\rm pr}_2^*{\cal G},{\rm pr}_1^!{\cal F})
\to R{\cal H}om(
{\rm pr}_2^*g_!{\cal G},{\rm pr}_1^!f_*{\cal F})
\label{eqf*}
\end{equation}
is defined in \cite{SGA5} (3.3.1),
using the isomorphism (\ref{eqFG}).
The diagram (\ref{eqpush})
defines a commutative diagram
\begin{equation}
\begin{CD}
C@>{c}>> X\times Y\\
@VhVV @VV{f\times g}V\\
C'@>{c'}>> X'\times Y'.
\end{CD}
\end{equation}

We assume the vertical arrows 
in (\ref{eqpush}) are proper.
Let 
$u\colon c_2^*{\cal G}\to c_1^!{\cal F}$
be a cohomological correspondence.
We identify $u$ with a 
section of $c^!R{\cal H}om(
{\rm pr}_2^*{\cal G},{\rm pr}_1^!{\cal F})$.
Then, it defines a section $h_*u$ of
$h_*c^!R{\cal H}om(
{\rm pr}_2^*{\cal G},{\rm pr}_1^!{\cal F})$.
By the assumption that
$f,g$ and $h$ are proper,
the base change map (\ref{f^!g_!})
defines a map of functors
$h_*c^!=
h_!c^!\to 
c^{\prime !}(f\times g)_!=
c^{\prime !}(f\times g)_*$.
This and the isomorphism (\ref{eqf*})
give maps
\begin{eqnarray*}
h_*c^!R{\cal H}om(
{\rm pr}_2^*{\cal G},{\rm pr}_1^!{\cal F})
&\to& c^{\prime !}(f\times g)_*
R{\cal H}om(
{\rm pr}_2^*{\cal G},{\rm pr}_1^!{\cal F})\\
&\to& c^{\prime !}R{\cal H}om(
{\rm pr}_2^*g_!{\cal G},{\rm pr}_1^!f_*{\cal F}).
\end{eqnarray*}
The push-forward
$h_*u\colon 
c_2^{\prime *}g_*{\cal G}=
c_2^{\prime *}g_!{\cal G}
\to c_1^{\prime !}f_*{\cal F}$
is defined
by the image of $h_*u$
by the composition.
The push-forward $h_*u$
is equal to the composition of the maps
$$\begin{CD}
c^{\prime*}_2g_*{\cal G}
@>{(\ref{g^*f_*})}>>
 h_*c^*_2{\cal G}
@>{h_*(u)}>> 
h_*c^!_1{\cal F}
@>{(\ref{f^!g_!})}>> 
c^{\prime!}_1f_*{\cal F}.
\end{CD}$$
The formation of the push-forward
is compatible with the composition.

A variant $h_!u$ of $h_*u$ 
is defined if the map
$c_2\colon C\to V$ in
(\ref{eqpush}) is proper.
In this case, we define $h_!u$ to be the composition
$$\begin{CD}
c^{\prime*}_2g_!{\cal G}
@>{(\ref{g^*f_!})}>>
 h_!c^*_2{\cal G}
@>{h_!(u)}>> 
h_!c^!_1{\cal F}
@>{(\ref{f^!g_!})}>> 
c^{\prime!}_1f_!{\cal F}.
\end{CD}$$
The formation of the $h_!u$
is also compatible with the composition.
If $f,g,h$ and $c_2$ are proper,
we have $h_*u=h_!u$.

We study the restriction
of a cohomological correspondence
to an open subscheme.
We consider a commutative diagram
\begin{equation}
\begin{CD}
U@<{c_1}<< C@>{c_2}>> V\\
@V{j_U}VV @VV{j_C}V @VV{j_V}V\\
X@<{\bar c_1}<< \overline C@>{\bar c_2}>> Y,
\end{CD}
\label{eqjU}
\end{equation}
of schemes over $k$
where the vertical arrows are 
open immersions.
Let ${\cal F}$ and ${\cal G}$
be objects of
$D_{\rm ctf}(X)$ and of
$D_{\rm ctf}(Y)$ respectively and
$\bar u\colon \bar c_2^*{\cal G}\to \bar c_1^!{\cal F}$
be a cohomological correspondence on $\overline C$.
Let ${\cal F}_U=j_U^*{\cal F}$
and ${\cal G}_V=j_V^*{\cal G}$
be the restrictions.
We identify
$j_C^*\bar c_1^!{\cal F}=
j_C^!\bar c_1^!{\cal F}=
c_1^!{\cal F}_U$
by the composition isomorphism.
Then, the restriction
$j_C^*\bar u$ on $C$
defines a cohomological correspondence
$u\colon c_2^*{\cal G}_V=
j_C^*\bar c_2^*{\cal G}
\to j_C^*\bar c_1^!{\cal F}=
c_1^!{\cal F}_U$.

\begin{lm}\label{lmopen}
Let the notation be as above
and let $j\colon U\times V\to X\times Y$
be the product $j_U\times j_V$.
We put ${\cal H}=
R{\cal H}om(
{\rm pr}_2^*{\cal G}_V,{\rm pr}_1^!{\cal F}_U)$
on $U\times V$
and $\overline{\cal H}=
R{\cal H}om(
{\rm pr}_2^*{\cal G},{\rm pr}_1^!{\cal F})$
on $X\times Y$.
We identify a cohomological correspondence
$\bar u\colon \bar c_2^*{\cal G}\to \bar c_1^!{\cal F}$
with a section
$\bar u$ of $\bar c^!\overline{\cal H}$
and the associated map
$\bar u\colon \bar c_!\Lambda \to \overline{\cal H}$.
We also identify the restriction
$u=j_C^*{\bar u}\colon 
c_2^*{\cal G}_V\to c_1^!{\cal F}_U$
with a section 
$u$ of $c^!{\cal H}$
and the associated map
$u\colon c_!\Lambda \to {\cal H}$.
Then, we have the following.

1. The section
$u$ of $c^!{\cal H}$
is the image of the restriction of $\bar u$
by the composition isomorphism
$j_C^*\bar c^!\overline{\cal H}=
j_C^!\bar c^!\overline{\cal H}
\to 
c^!j^!\overline{\cal H}
=c^!j^*\overline{\cal H}
=c^!{\cal H}.$

2. The square
\begin{equation}
\begin{CD}
j^*\bar c_!\Lambda @>{j^*\bar u}>> 
j^*\overline{\cal H}\\
@A{\rm (\ref{f^!g_!})}AA @VV{(\ref{eqw!})}V\\
c_!\Lambda @>u>>{\cal H}
\end{CD}\label{eqopen}
\end{equation}
is commutative.
\end{lm}

{\it Proof.}
1. The isomorphisms (\ref{eqw!})
form a commutative diagram
$$\begin{CD}
j_C^!\bar c^!R{\cal H}om(
{\rm pr}_2^*{\cal G},{\rm pr}_1^!{\cal F})
@>>> 
j_C^!R{\cal H}om(
\bar c_2^*{\cal G},\bar c_1^!{\cal F})\\
@VVV @VVV\\
c^!R{\cal H}om(
{\rm pr}_2^*{\cal G}_V,{\rm pr}_1^!{\cal F}_U)
@>>> 
R{\cal H}om(
c_2^*{\cal G}_V,\bar c_1^!{\cal F}_U).
\end{CD}
$$
The assertion follows from this.

2. 
We consider the diagram
$$
\begin{CD}
j^*\bar c_!\Lambda 
@>{j^*\bar c_!\bar u}>> 
j^*\bar c_!\bar c^!\overline{\cal H}
@>{\rm adj}>> 
j^*\overline{\cal H}\\
@A{\rm (\ref{f^!g_!})}AA 
@AA{\rm (\ref{f^!g_!})}A\\
c_!j_C^*\Lambda 
@>{c_!j_C^*\bar u}>> 
c_!j_C^*\bar c^!\overline{\cal H}
@. @VV{(\ref{eqw!})}V\\
@V{\simeq}VV @V{\simeq}VV @.\\
c_!\Lambda 
@>{c_!u}>> c_!c^!{\cal H}
@>{\rm adj}>>{\cal H}.
\end{CD}$$
The upper left square is
commutative by functoriality.
The lower left square is
commutative by 1
and the vertical arrows are
isomorphisms.
The right rectangle is
commutative by the definition
of the base change map (\ref{f^!g_!}).
Since the square (\ref{eqopen}) is
the outer square,
the assertion is proved.
\qed

\begin{lm}\label{lmext}
Assume that the right square 
in the diagram {\rm (\ref{eqjU})}
is cartesian.
Let ${\cal F}$ and ${\cal G}$
be objects of
$D_{\rm ctf}(U)$ and of
$D_{\rm ctf}(V)$ respectively and
$u\colon c_2^*{\cal G}\to c_1^!{\cal F}$
be a cohomological correspondence on $C$.
Then, there exists a unique cohomological
correspondence
$\bar u\colon \bar c^*_2j_{V!}{\cal G}
\to \bar c^!_1j_{U!}{\cal F}$
on $\overline C$
such that $j_C^*\bar u=u$.
\end{lm}

{\it Proof.}
Since the right square in (\ref{eqjU})
is assumed cartesian,
we have
$\bar c^*_2j_{V!}
=j_{C!}c_2^*$.
By adjunction,
there exists a unique map
$\bar c^*_2j_{V!}{\cal G}
=j_{C!}c_2^*{\cal G}
\to \bar c^!_1j_{U!}{\cal F}$
corresponding to
$c_2^*{\cal G}
\to j_C^*\bar c^!_1j_{U!}{\cal F}
=c^!_1j_C^*j_{U!}{\cal F}$.
\qed

\begin{cor}\label{corext}
1. Assume that
the map $c_2\colon C\to V$ is proper and
$C$ is dense in $\overline C$.
Then the right square 
in the diagram {\rm (\ref{eqjU})}
is cartesian.

2. Assume that the right square 
in the diagram {\rm (\ref{eqjU})}
is cartesian.
Let ${\cal F}$ and ${\cal G}$
be objects of
$D_{\rm ctf}(U)$ and of
$D_{\rm ctf}(V)$ respectively and
$u\colon c_2^*{\cal G}\to c_1^!{\cal F}$
be a cohomological correspondence on $C$.
Then, $\bar u=j_{C!}u\colon \bar c^*_2j_{V!}{\cal G}
\to \bar c^!_1j_{U!}{\cal F}$ is the 
unique cohomological correspondence
on $\overline C$
such that $j_C^*\bar u=u$.
\end{cor}

{\it Proof.}
1. Since $C$ is closed and dense in
$\overline C\times_YV$,
we have $C=\overline C\times_YV$.

2. Clear from $j_C^*j_{C!}u=u$
and the uniqueness proved in Lemma \ref{lmext}.
\qed

We call $j_{C!}u\colon \bar c^*_2j_{V!}{\cal G}
\to \bar c^!_1j_{U!}{\cal F}$
the zero-extension of $u$
(cf. \cite{Pink} (2.3)).
Similarly as the isomorphism
(\ref{eqf*}),
the isomorphism (\ref{eqFG})
induces a canonical isomorphism
\begin{equation}
(1\times j_V)_*
(j_U\times 1)_!
R{\cal H}om({\rm pr}_2^*{\cal G},
{\rm pr}_1^!{\cal F})\to
R{\cal H}om({\rm pr}_2^*j_{V!}{\cal G},
{\rm pr}_1^!j_{U!}{\cal F}).
\label{eqHj}
\end{equation}
If the map
$\bar c\colon \overline C\to X\times Y$
is a closed immersion and
if $u\colon c^*_2{\cal G}
\to c^!_1{\cal F}$
is identified with
an element $u\in H^0_C(U\times V,
R{\cal H}om({\rm pr}^*_2{\cal G},
{\rm pr}^!_1{\cal F}))$,
then the zero-extension
$j_{C!}u$ is the 
inverse image of $u$ 
by the isomorphism
\begin{equation}
\begin{array}{cl}
&H^0_{\overline C}(X\times Y,
R{\cal H}om({\rm pr}^*_2j_{V!}{\cal G},
{\rm pr}^!_1j_{U!}{\cal F}))\\
&\\
\overset \sim\gets&
H^0_{\overline C}(X\times Y,
(1\times j_V)_*
(j_U\times 1)_!
R{\cal H}om({\rm pr}_2^*{\cal G},
{\rm pr}_1^!{\cal F}))\\
&\\
\to&
H^0_C(X\times V,
(j_U\times 1)_!
R{\cal H}om({\rm pr}_2^*{\cal G},
{\rm pr}_1^!{\cal F}))
=
H^0_C(U\times V,
R{\cal H}om({\rm pr}^*_2{\cal G},
{\rm pr}^!_1{\cal F})).
\end{array}
\label{eqzeroex}
\end{equation}

We define the pull-back
of a cohomological correspondence.
Let $f\colon X'\to X$ and $g\colon Y'\to Y$
be morphisms of smooth schemes 
over $k$.
We assume 
$\dim X=\dim X'$ and $\dim Y=\dim Y'$.
Let ${\cal F}$ and ${\cal G}$
be objects of
$D_{\rm ctf}(X)$
and of
$D_{\rm ctf}(Y)$
respectively.
Then the canonical maps
(\ref{f^*f^!}) and 
(\ref{eqw!})
induce a map
\begin{equation}
\begin{CD}
(f\times g)^*R{\cal H}om(
{\rm pr}_2^*{\cal G},{\rm pr}_1^!{\cal F})
\to &
(f\times g)^!R{\cal H}om(
{\rm pr}_2^*{\cal G},{\rm pr}_1^!{\cal F})\\
&\to 
R{\cal H}om(
{\rm pr}_2^{\prime *}g^*{\cal G},{\rm pr}_1^{\prime !}f^!{\cal F}).
\end{CD}
\label{eqf^*}
\end{equation}
With the isomorphism
(\ref{eqFG}),
the map (\ref{eqf^*})
is identified with the composition
\begin{equation}
\begin{CD}
(f\times g)^*
({\cal F}\boxtimes {\mathbf D}{\cal G})
\to 
f^*{\cal F}\boxtimes g^*{\mathbf D}{\cal G}
\to 
f^!{\cal F}\boxtimes g^!{\mathbf D}{\cal G}
\to 
f^!{\cal F}\boxtimes {\mathbf D}g^*{\cal G}
\end{CD}
\label{eqf^**}
\end{equation}
where the middle arrow is defined by 
(\ref{f^*f^!}).

Let 
$c=(c_1,c_2)\colon C\to X\times Y$
be a correspondence
and 
$u\colon c_2^*{\cal G}\to c_1^!{\cal F}$
be a cohomological correspondence 
on $C$.
We identify $u$ with a map
$u\colon c_!\Lambda\to 
R{\cal H}om(
{\rm pr}_2^*{\cal G},{\rm pr}_1^!{\cal F})$ as above.
We define 
a correspondence
$c'=(c'_1,c'_2)\colon C'\to X'\times Y'$
by the cartesian diagram
\begin{equation}
\begin{CD}
C'@>{c'}>> X'\times Y'\\
@VhVV @VV{f\times g}V\\
C@>{c}>> X\times Y.
\end{CD}
\end{equation}
By the proper base change theorem, 
the base change map
$(f\times g)^*c_!\Lambda\to c'_!\Lambda$
is an isomorphism.
Hence the map
$u\colon c_!\Lambda\to R{\cal H}om(
{\rm pr}_2^*{\cal G},{\rm pr}_1^!{\cal F})$ induces 
a map
$$\begin{CD}
c'_!\Lambda \overset\sim\gets (f\times g)^*c_!\Lambda
\to
(f\times g)^*R{\cal H}om(
{\rm pr}_2^*{\cal G},{\rm pr}_1^!{\cal F})
@>{(\ref{eqf^*})}>>
R{\cal H}om(
{\rm pr}_2^{\prime *}g^*{\cal G},{\rm pr}_1^{\prime !}f^!{\cal F}).
\end{CD}$$
The composition defines a cohomological correspondence
$(f\times g)^*u\colon c_2^{\prime *}g^*{\cal G}
\to c_1^{\prime !}f^!{\cal F}=
c_1^{\prime !}f^*{\cal F}$.
We call $(f\times g)^*u$ the pull-back 
of $u$ by $f\times g$.
If $c$ is a closed immersion,
the correspondence $u\colon c_2^*{\cal G}
\to c_1^!{\cal F}$
is identified with a cohomology class
$u\in H^0_C(X\times Y,
R{\cal H}om(
{\rm pr}_2^*{\cal G},{\rm pr}_1^!{\cal F}))$
and the pull-back $(f\times g)^*u$
is identified
with the pull-back
$(f\times g)^*u\in H^0_{C'}(X'\times Y',
R{\cal H}om(
{\rm pr}_2^{\prime*}g^*{\cal G},
{\rm pr}_1^{\prime!}f^!{\cal F}))$.

For a cohomological correspondence
$u(\Gamma,\gamma)$,
its pull-back is
computed by using the intersection product.

\begin{lm}\label{lmupull}
Let the notation be as above.
The diagram
\begin{equation}
\begin{CD}
CH_d(C)\otimes
Hom(c_2^*{\cal G},c_1^*{\cal F})
@>{\rm(\ref{equgam})}>> H^0_C(X\times Y,
R{\cal H}om({\rm pr}_2^*{\cal G},
{\rm pr}_1^!{\cal F}))\\
@V{(f\times g)^!\otimes h^*}VV 
@VV{(f\times g)^*}V\\
CH_d(C')\otimes
Hom(c_2^{\prime *}g^*{\cal G},c_1^{\prime *}f^*{\cal F})
@>{\rm(\ref{equgam})}>> 
H^0_{C'}(X'\times Y',
R{\cal H}om({\rm pr}_2^*g^*{\cal G},
{\rm pr}_1^!f^*{\cal F}))
\end{CD}
\label{equpull}
\end{equation}
is commutative.
\end{lm}

{\it Proof.}
Clear from the commutative diagram
(\ref{eqcycpull}).
\qed

\subsection{Complement on the Verdier pairing.}

We briefly recall the definition 
of the pairing \cite{SGA5} (4.2.5).
We consider a cartesian diagram
$$
\xymatrix{
D\ar[d]_d& E\ar[l]_{c'}\ar[d]^{d'}\ar[dl]_e\\
X&C\ar[l]_c}$$
of schemes over $k$.
Let ${\cal P}$ and ${\cal Q}$
be objects of $D_{\rm ctf}(X)$
and let ${\cal P}\otimes {\cal Q}
\to {\cal K}_X$ be a morphism.
Then, the map (4.2.1) in \cite{SGA5} and
the map ${\cal P}\otimes {\cal Q}
\to {\cal K}_X$
induce a map
$c^!{\cal P}\boxtimes d^!{\cal Q}
\to e^!{\cal K}_X={\cal K}_E$
and hence a pairing
\begin{equation}
\langle\ ,\ \rangle\colon 
H^0(C,c^!{\cal P})\otimes 
H^0(D,d^!{\cal Q})
\to H^0(E,e^!{\cal K}_X)=H^0(E,{\cal K}_E).
\label{eqVer}
\end{equation}
If the maps
$c\colon C\to X$ and
$d\colon D\to X$ are closed immersions,
the pairing (\ref{eqVer}) gives a 
pairing
$\langle\ ,\ \rangle\colon 
H^0_C(X,{\cal P})\otimes 
H^0_D(X,{\cal Q})
\to H^0_E(X,{\cal K}_X)=H^0(E,{\cal K}_E)$.

\begin{lm}\label{lmcomp}
Let the notation be as above.
Then, we have a commutative diagram
$$\begin{CD}
H^0(C,c^!{\cal P})\otimes 
H^0(D,d^!{\cal Q})
@>{\langle\ ,\ \rangle}>>
H^0(E,{\cal K}_E)\\
@V{d^{\prime*}\otimes 1}VV @|\\
H^0(E,c^{\prime !}d^*{\cal P})\otimes 
H^0(D,d^!{\cal Q})
@>{\langle\ ,\ \rangle}>>
H^0(E,{\cal K}_E)
\end{CD}$$
where the lower horizontal arrow
is the pairing 
{\rm (\ref{eqVer})} for 
$E\overset{c'}\to D
\overset{\rm id}\gets D$ and
the map
$d^*{\cal P}
\otimes
d^!{\cal Q}
\to 
d^!{\cal K}_X={\cal K}_D$
{\rm (\ref{FoxG})}
induced by
${\cal P}\otimes
{\cal Q}\to {\cal K}_X$.
\end{lm}

The pairing $\langle\ ,\ \rangle$
is compatible with the proper push-forward in
the following sense. 

\begin{lm}\label{lmVer}
Let 
\begin{equation}
\xymatrix{
 &D\ar[dl]_{d}\ar[dd]_-(0.7){f_D}&&E\ar[dl]\ar[ll]\ar[dd]\\
X\ar[dd]_f&&C\ar[ll]_{\qquad c}\ar[dd]_-(0.7){f_C}&\\
 &D'\ar[dl]_{d'}&&E'\ar[dl]\ar[ll]\\
X'&&C'\ar[ll]_{\quad c'}&}
\end{equation}
be a commutative cube of
schemes over $k$.
We assume that the horizontal faces
are cartesian
and that vertical arrows are proper.
Let ${\cal P}$
and ${\cal Q}$
be objects of $D_{\rm ctf}(X)$
and let ${\cal P}\otimes {\cal Q}
\to {\cal K}_X$ be a morphism.

Then, the pairings $\langle\ ,\ \rangle$
for ${\cal P}\otimes {\cal Q}
\to {\cal K}_X$
and for 
$f_*{\cal P}\otimes f_*{\cal Q}
\to f_*{\cal K}_X\to {\cal K}_{X'}$
form a commutative diagram
\begin{equation}
\begin{CD}
H^0(C,c^!{\cal P})\otimes 
H^0(D,d^!{\cal Q})
@>{\langle\ ,\ \rangle}>> H^0(E,{\cal K}_E)\\
@VVV @VVV\\
H^0(C',c^{\prime !}f_*{\cal P})\otimes 
H^0(D',d^{\prime !}f_*{\cal Q})
@>{\langle\ ,\ \rangle}>> H^0(E',{\cal K}_{E'})\\
\end{CD}
\label{eqVp}
\end{equation}
where the left vertical maps
are induced by the base change maps
$f_{C*}c^!\to c^{\prime !}f_*$
and
$f_{D*}d^!\to d^{\prime !}f_*$
{\rm (\ref{f^!g_!})}.
\end{lm}

{\it Proof.}
It follows from the commutativity of
the squares (A) and (B)
in the diagram \cite{SGA5} (4.4.2).
\qed

Applying the pairing (\ref{eqVer})
to cohomological correspondences,
we obtain the Verdier pairing.
Let $X$ and $Y$ be schemes
over $k$ and 
$$\begin{CD}
C@<<< E\\
@V{c=(c_1,c_2)}VV @VVV\\
X\times Y@<{d=(d_1,d_2)}<<D
\end{CD}$$
be a cartesian diagram
of schemes over $k$.
Let ${\cal F}$ and ${\cal G}$
be objects of
$D_{\rm ctf}(X)$ and of
$D_{\rm ctf}(Y)$ respectively.
We put 
${\cal H}=R{\cal H}om(
{\rm pr}_2^*{\cal G},{\rm pr}_1^!{\cal F})$
and 
${\cal H}^*=R{\cal H}om(
{\rm pr}_1^*{\cal F},{\rm pr}_2^!{\cal G})$
on $X\times Y$.
A canonical map
\begin{equation}
\begin{CD}
{\cal H}\otimes {\cal H}^*\to
{\cal K}_{X\times Y}
\end{CD}
\label{eqVerdi}
\end{equation}
is defined in \cite{SGA5} (4.1.4),
using the isomorphism (\ref{eqFG}).
If $X$ and $Y$ are smooth over $k$
of dimension $d$
and if ${\cal F}$ and
${\cal G}$ are smooth sheaves
of free $\Lambda$-modules,
we have
${\cal H}={\cal H}om(
{\rm pr}_2^*{\cal G},{\rm pr}_1^*{\cal F})
(d)[2d]$
and 
${\cal H}^*={\cal H}om(
{\rm pr}_1^*{\cal F},{\rm pr}_2^*{\cal G})
(d)[2d]$
and the map
(\ref{eqVerdi})
is induced by the pairing
${\cal H}om(
{\rm pr}_2^*{\cal G},{\rm pr}_1^*{\cal F})
\otimes {\cal H}om(
{\rm pr}_1^*{\cal F},{\rm pr}_2^*{\cal G})
\to \Lambda\colon 
a\otimes b\mapsto {\rm Tr}(b\circ a)$.
The pairing (\ref{eqVer})
defines the Verdier pairing
\begin{equation}
\langle\ ,\ \rangle\colon 
Hom(c_2^*{\cal G},c_1^!{\cal F})
\otimes 
Hom(d_2^*{\cal F},d_1^!{\cal G})
\to H^0(E,{\cal K}_E)
\label{eqVerd}
\end{equation}
as in \cite{SGA5} (4.2.5).

The following result will be useful
in computing the Verdier pairing
of the extension by zero of
a correspondence.

\begin{pr}\label{lmpink}
Let $X$ and $Y$ be schemes over $k$
and $j_U\colon U\to X$ 
and $j_V\colon V\to Y$ be open immersions.
Let $f\colon (X\times Y)'\to X\times Y$
be a proper morphism over $k$
that is an isomorphism on 
$(X\times V)\cup (U\times Y)$;
we identify
$(X\times V)\cup (U\times Y)$
as an open subscheme of
$(X\times Y)'$ by the restriction of $f$.
We consider a commutative diagram
\begin{equation}
\xymatrix{
& C\ar[dl]_{j_C'}\ar[dddl]^{j_C}\ar[rr]^{\!\!\!\! c}
&&U\times V\ar[dl]_{j'}\ar[dddl]^{j=j_U\times j_V}&&
D\ar[dl]_{j_D'}\ar[dddl]^{j_D}\ar[ll]_{\quad d}\\
{\overline C}'\ar[dd]_{g}\ar[rr]^{\!\!\!\! \bar c'}&&
(X\times Y)'
\ar[dd]_{f}&&
{\overline D}'\ar[dd]_{h}\ar[ll]_{\qquad \bar d'}
&\\\\
{\overline C}\ar[rr]^{\!\!\!\! \bar c}&&
X\times Y&&
{\overline D}\ar[ll]_{\quad \bar d}&}
\label{eqpink}
\end{equation}
of schemes over $k$.
We assume that the vertial arrows $f,g,h$
and the horizontal arrows $c,\bar c,\bar c',
d,\bar d,\bar d'$
are proper
and that the four slant parallelograms are cartesian,
so that the six slant arrows
are open immersions.
We also assume
$C=
\overline C\times_{X\times Y}(X\times V)$
and 
$D=
\overline D\times_{X\times Y}(U\times Y)$.

Let ${\cal F}$ be an object of
$D_{\rm ctf}(U)$
and ${\cal G}$ be an object of
$D_{\rm ctf}(V)$ respectively.
We put 
${\cal H}=R{\cal H}om(
{\rm pr}_2^*{\cal G},{\rm pr}_1^!{\cal F})$
and 
${\cal H}^*=R{\cal H}om(
{\rm pr}_1^*{\cal F},{\rm pr}_2^!{\cal G})$
on $U\times V$
and 
$\overline {\cal H}=R{\cal H}om(
{\rm pr}_2^*j_{V!}{\cal G},{\rm pr}_1^!j_{U!}{\cal F})$
and 
$\overline {\cal H}^*=R{\cal H}om(
{\rm pr}_1^*j_{U!}{\cal F},{\rm pr}_2^!j_{V!}{\cal G})$
on $X\times Y$.
Let
$\overline {\cal H}'$
and 
$\overline {\cal H}^{*\prime}$
be objects of
$D_{\rm ctf}((X\times Y)')$.
Let
$u\colon c_!\Lambda\to {\cal H}$
and $v\colon d_!\Lambda\to {\cal H}^*$
be cohomological correspondences.

1.
For morphisms
$a\colon \overline {\cal H}
|_{U\times Y}
\to
\overline {\cal H}'
|_{U\times Y}$
and 
$b\colon \overline {\cal H}^*
|_{X\times V}
\to
\overline {\cal H}^{*\prime}
|_{X\times V}$,
there exist unique morphisms
$\bar a\colon \overline {\cal H}
\to
f_*\overline {\cal H}'$
and 
$\bar b\colon \overline {\cal H}^*
\to
f_*\overline {\cal H}^{*\prime}$
extending $a$ and $b$ respectively.

2. 
Let $\bar u'\colon 
\bar c'_!\Lambda \to
\overline {\cal H}'$
and $\bar v'\colon 
\bar d'_!\Lambda \to
\overline {\cal H}^{\prime*}$
be morphisms
and $\bar a\colon \overline {\cal H}
\to f_*\overline {\cal H}'$
and $\bar b\colon \overline {\cal H}^*
\to f_*\overline {\cal H}^{*\prime}$
be isomorphisms.
We identify the functor
$j^*f_*$ with $j^{\prime *}$
by the isomorphism 
$j^{\prime *}\to j^*f_*$
{\rm (\ref{f^!g_!})}.
Suppose the squares
\begin{equation}
\begin{CD}
c_!\Lambda @>u>> {\cal H}\\
@V{\rm (\ref{f^!g_!})}VV 
@VV{j^*\bar a}V\\
j^{\prime*}\bar c'_!\Lambda 
@>{j^{\prime*}\bar u'}>> 
j^{\prime*}\overline{\cal H}'
\end{CD}\qquad\qquad
\begin{CD}
d_!\Lambda @>v>> {\cal H}^*\\
@V{\rm (\ref{f^!g_!})}VV 
@VV{j^*\bar b}V\\
j^{\prime*}\bar d'_!\Lambda 
@>{j^{\prime*}\bar v'}>> 
j^{\prime*}\overline{\cal H}^{*\prime}
\end{CD}
\label{eqlmpink1}\end{equation}
are commutative.
Then, for the zero-extensions
$j_{C!}u$ and $j_{D!}v$, the squares
\begin{equation}
\begin{CD}
\bar c_!\Lambda @>{j_{C!}u}>> \overline{\cal H}\\
@VVV @VV{\bar a}V\\
f_*\bar c'_!\Lambda @>{f_*\bar u'}>> 
f_*\overline{\cal H}'
\end{CD}\qquad\qquad
\begin{CD}
\bar d_!\Lambda @>{j_{D!}v}>> \overline{\cal H}^*\\
@VVV @VV{\bar b}V\\
f_*\bar d'_!\Lambda 
@>{f_*\bar v'}>> 
f_*\overline{\cal H}^{*\prime}
\end{CD}
\label{eqlmpink2}\end{equation}
are commutative,
where the left vertical arrows
are induced by the canonical maps
$\Lambda\to g_*\Lambda$ and
$\Lambda\to h_*\Lambda$.

3. 
We keep the assumpions in 2.
Further, let
$\overline {\cal H}'
\otimes 
\overline {\cal H}^{*\prime}
\to {\cal K}_{(X\times Y)'}$
be a map
making the square
\begin{equation}
\begin{CD}
{\cal H}
\otimes 
{\cal H}^*
@>{\rm (\ref{eqVerdi})}>> {\cal K}_{U\times V}\\
@V{j^*\bar a\otimes j^*\bar b}VV
@AA{\rm can.}A\\
j^{\prime*}
(\overline {\cal H}'
\otimes 
\overline {\cal H}^{*\prime})
@>>> j^{\prime*}{\cal K}_{(X\times Y)'}
\end{CD}
\label{eqpinkp}
\end{equation}
commutative.
Then, 
we have
$$\langle j_{C!}u,j_{D!}v\rangle=
f_*\langle \bar u',\bar v'\rangle$$
in $H^0(E,{\cal K}_E)$.
\end{pr}

{\it Proof.}
1.
Since
$\overline {\cal H}
=(j_U\times 1)_!{\cal H}|_{U\times Y}$
and 
$\overline {\cal H}^*
=(1\times j_V)_!{\cal H}^*|_{X\times V}$,
the assertion follows.

2.
We show the assertion for $u$.
The proof for $v$ is similar and omitted.
We claim that we obtain
(\ref{eqlmpink1})
by applying $j^*$
to (\ref{eqlmpink2}).
This is clear for the bottom line and the right column.
For the top line,
we have $j_C^*j_{C!}u=u$
by Corollary \ref{corext}.2.
We show it for the left column.
It suffices to show the diagram
$$\xymatrix{
j^*\bar c_!\Lambda 
\ar[d]&
c_!\Lambda 
\ar[l]_{(\ref{f^!g_!})}
\ar[dl]_{(\ref{f^!g_!})}
\ar[d]^{(\ref{f^!g_!})}\\
j^*f_*\bar c'_!\Lambda &
\ar[l]_{(\ref{f^!g_!})}
j^{\prime*}\bar c'_!\Lambda 
}$$
is commutative.
The lower right triangle
is commutative by the transitivity
(\ref{f^!g_!2}) of the base change maps
(\ref{f^!g_!})
for the base change of
$\overline C\to (X\times Y)'
\to X\times Y$
by $j\colon U\times V\to X\times Y$.
Since the base change map
(\ref{f^!g_!})
is the inverse of
(\ref{g^*f_*}),
the upper left triangle
is the adjoint of
the transitivity for 
the canonical maps
$\bar c_!\Lambda 
\to
f_*\bar c'_!\Lambda 
\to
j_*c_!\Lambda$
and is commutative.
Hence,
the left column in
(\ref{eqlmpink2}) also induces
that in
(\ref{eqlmpink1}).

Since $\bar a$ is an isomorphism,
there exists a unique map
$\bar u\colon \bar c_!\Lambda
\to \overline{\cal H}$
that makes the square
(\ref{eqlmpink2})
commutative.
Since $u$ is the unique map
that makes the square
(\ref{eqlmpink1})
commutative,
we have 
$j_C^*\bar u=u$.
Thus we obtain $\bar u=j_{C!}u$
by Lemma \ref{lmext}.

3.
We defined the push forward 
$f_*\bar u'\colon 
\Lambda\to
g_*\bar c^{\prime!}\overline{\cal H}'\to
\bar c^!f_*\overline{\cal H}'$
and proved
$f_*\langle \bar u',\bar v'\rangle=
\langle f_*\bar u',f_*\bar v'\rangle$
in Lemma \ref{lmVer}.
The push forward 
$f_*\bar u'$ is the adjoint of 
the composition
$\bar c_!\Lambda\to
f_*\bar c'_!\Lambda
\overset{f_*\bar u'}\to
f_*\overline{\cal H}'$
and is equal to
$\bar a\circ j_{C!}u$ by 2.
Hence, 
we have
$\langle f_*\bar u',f_*\bar v'\rangle=
\langle \bar a\circ j_{C!}u,
\bar b \circ j_{D!}v\rangle$.
Thus, it suffices to show
$\langle j_{C!}u,j_{D!}v\rangle
=
\langle \bar a\circ j_{C!}u,
\bar b \circ j_{D!}v\rangle$.

Since $\overline {\cal H}
\otimes 
\overline {\cal H}^*
=j_!
({\cal H}\otimes 
{\cal H}^*)$,
we obtain a commutative diagram
\begin{equation}
\begin{CD}
\overline {\cal H}
\otimes 
\overline {\cal H}^*
@>{(\ref{eqVerdi})}>>
{\cal K}_{X\times Y}\\
@V{\bar a\otimes \bar b}VV @AAA\\
f_*\overline {\cal H}'
\otimes 
f_*\overline {\cal H}^{*\prime}
@>>> f_*{\cal K}_{(X\times Y)'}
\end{CD}
\label{eqpv}
\end{equation}
where the lower horizontal arrow
is induced by
$\overline {\cal H}'
\otimes 
\overline {\cal H}^{*\prime}
\to {\cal K}_{(X\times Y)'}$.
Thus, 
we have
$\langle j_{C!}u,j_{D!}v\rangle
=
\langle \bar a\circ j_{C!}u,
\bar b \circ j_{D!}v\rangle$
as required.
\qed

\section{Characteristic class
of an $\ell$-adic \'etale sheaf.}

\subsection{Characteristic class
of a cohomological correpondence.}

Let $X$ be a scheme over $k$ and
$\delta\colon 
X=\Delta_X\to X\times X$ be the diagonal map.
Let ${\cal F}$ be an object of $D_{\rm ctf}(X)$
and let $1=u(X,1)$ be
the cohomological correspondence
defined by the identity of ${\cal F}$ on 
the diagonal $X$.

\begin{df}\label{dfcc}
Let $c=(c_1,c_2)\colon C\to X\times X$
be a closed immersion and
$u\colon c_2^*{\cal F}\to c_1^!{\cal F}$
be a cohomological correspondence on $C$.
Then, we call the cohomology class
$\langle u,1\rangle
\in H^0_{C\cap X}(X,{\cal K}_X)$
defined in {\rm (\ref{eqVerd})}
{\rm the characteristic class}
and write
$C({\cal F},C,u)$.

If $C=X$ is the diagonal
and $u\colon {\cal F}\to {\cal F}$
is an endomorphism,
we drop $C$ from the notation and we
write $C({\cal F},u)\in H^0(X,{\cal K}_X)$.
Further, if $u$ is the identity,
we simply write 
$C({\cal F})=\langle 1,1\rangle$
and call it the characteristic class of ${\cal F}$.
\end{df}

\begin{eg}[\cite{ksc} Exercise I.32]
If $X=C={\rm Spec}\ k$,
we have
$C({\cal F},u)={\rm Tr}\ u
\in \Lambda$.
In particular,
if $u=1$,
we have
$C({\cal F})=\chi({\cal F}).$
\end{eg}

\begin{lm}\label{lmadd}
Let $X$ be a scheme over $k$
and $c=(c_1,c_2)\colon C\to X\times X$
be a closed immersion.
Let $({\cal F},F_\bullet)$
be an object of the
filtered derived category
$DF_{\rm ctf}(X)$
and 
$u\colon c_2^*{\cal F}\to
c_1^!{\cal F}$
be a morphism in
the filtered derived category
$DF_{\rm ctf}(C)$.
Then, we have
$$C({\cal F},C,u)=
\sum_qC({\rm Gr}^F_q {\cal F},C,
{\rm Gr}^F_qu)$$
in $H^0(X,{\cal K}_X)$.
\end{lm}

{\it Proof.}
It is a special case of 
the equality \cite{SGA5} (4.13.1).
For later use,
we briefly skecth the proof.
By the isomorphism (\ref{fHomg}),
the homomorphism
$u\in Hom_{DF_{\rm ctf}(C)}(
c_2^*{\cal F},c_1^!{\cal F})$
defines a section
$u$ of $c^!F_0R{\cal H}om(
{\rm pr}_2^*{\cal F},
{\rm pr}_1^!{\cal F})$ on $C$.
Similarly the identity
defines a section
$1$ of $\delta^!F_0R{\cal H}om(
{\rm pr}_1^*{\cal F},
{\rm pr}_2^!{\cal F})$ on 
the diagonal $X$.
We have a canonical isomorphism
$Gr^F_0R{\cal H}om(
{\rm pr}_2^*{\cal F},
{\rm pr}_1^!{\cal F})
\to
\bigoplus_q R{\cal H}om(
{\rm pr}_2^*Gr^F_q{\cal F},
{\rm pr}_1^!Gr^F_q{\cal F})$
and similarly for
$Gr^F_0R{\cal H}om(
{\rm pr}_1^*{\cal F},
{\rm pr}_2^!{\cal F})$.
Since the restriction
of the canonical pairing
(\ref{eqVer}) on
$F_0R{\cal H}om(
{\rm pr}_2^*{\cal F},
{\rm pr}_1^!{\cal F})
\otimes F_0R{\cal H}om(
{\rm pr}_1^*{\cal F},
{\rm pr}_2^!{\cal F})$ 
is induced by the sum
of the canonical pairings
$R{\cal H}om(
{\rm pr}_2^*Gr^F_q{\cal F},
{\rm pr}_1^!Gr^F_q{\cal F})
\otimes R{\cal H}om(
{\rm pr}_1^*Gr^F_q{\cal F},
{\rm pr}_2^!Gr^F_q{\cal F})
\to {\cal K}_{X\times X}$,
the assertion follows.
\qed

\begin{rmk}[\cite{fer}]
For a distinguished triangle
$\to {\cal F}\overset u\to 
{\cal G}\overset v\to 
{\cal H}\overset w\to $
of perfect complexes of $\Lambda$-modules
and endomorphisms $f,g,h$
of ${\cal F},{\cal G},{\cal H}$
compatible with $u,v,w$,
there is a counterexample 
to the equality
${\rm Tr}(g)={\rm Tr}(f)+{\rm Tr}(h)$.
\end{rmk}

\begin{cor}\label{coradd}
Let $X$ be a scheme over $k$
and $({\cal F},F)$
be an object of 
$DF_{\rm ctf}(X)$.
Then, we have
$$C({\cal F})=\sum_qC(Gr^F_q{\cal F})$$
in $H^0(X,{\cal K}_X)$.
\end{cor}

\begin{pr}[\cite{SGA5} Th\'eor\`eme 4.4,
Corollaire 4.8]\label{prpush}
Let $f\colon X\to Y$ be a 
proper morphism of
schemes over $k$
and we consider a commutative diagram
$$\begin{CD}
C@>{(c_1,c_2)}>>X\times X\\
@VhVV @VV{f\times f}V\\
C'@>>> Y\times Y
\end{CD}$$
of schemes   over $k$
where the vertical arrows
are proper and the horizontal arrows are closed immersions.
Let $u\colon c_2^*{\cal F}\to c_1^!{\cal F}$ be
a cohomological correspondence.
Then we have
$$C(f_*{\cal F},C',h_*u)=f_*C({\cal F},C,u).$$
In particular,
if $X$ is proper over $k$,
we have
$${\rm Tr}(u^*\colon H^*(X_{\overline k},
{\cal F}))={\rm Tr}\
C({\cal F},C,u),$$
where ${\overline k}$
is a separable closure of $k$.
\end{pr}

By devissage using Lemma \ref{lmadd}
and Proposition \ref{prpush},
the computation of a characteric class
is reduced to the computation of
$C(j_!{\cal F})$ 
where $j\colon U\to X$
is an open immersion,
$U$ is smooth and
the cohomology sheaves
${\cal H}^q{\cal F}$ are locally constant on
$U$.
We will later compute
the characteristic class $C(j_!{\cal F})$
or rather the difference
$C(j_!{\cal F})-{\rm rank}\ {\cal F}\cdot
C(j_!\Lambda)$
in terms of the ramification of ${\cal F}$
along the boundary $X\setminus U$.

We give
an equivalent description of the characteristic class.
The canonical isomorphism (\ref{eqFG})
\begin{equation}
\begin{CD}
{\cal H}=R{\cal H}om(
{\rm pr}_2^*{\cal F},{\rm pr}_1^!{\cal F})
@>>>{\cal F}
\boxtimes^L
{\mathbf D}{\cal F}\end{CD}
\label{eqH}
\end{equation}
induces an isomorphism
\begin{equation}
\begin{CD}
\delta^*{\cal H}
@>>>{\cal F}
\otimes^L
{\mathbf D}{\cal F}.\end{CD}
\label{eqdH}
\end{equation}
Thus the evaluation map
\begin{equation}
\begin{CD}
{\cal F}
\otimes^L
{\mathbf D}{\cal F}
@>>>{\cal K}_X
\end{CD}
\label{eqev}
\end{equation}
\cite{SGA5} (2.2.2)
induces a map
\begin{equation}
\begin{CD}
\delta^*{\cal H}
@>>>{\cal K}_X.
\end{CD}
\label{eqeva}
\end{equation}
If $X$ is smooth of dimension $d$
and if ${\cal F}$ is a smooth sheaf
of free $\Lambda$-modules,
the map (\ref{eqeva}) is
equal to the composition of
$$\begin{CD}
\delta^*{\cal H}
=\delta^*{\cal H}om({\rm pr}_2^*{\cal F},
{\rm pr}_1^*{\cal F})(d)[2d]
={\cal E}nd({\cal F})(d)[2d]
@>{\rm Tr}>>
\Lambda(d)[2d]\to
{\cal K}_X.
\end{CD}$$
If we put
${\cal H}^*=R{\cal H}om(
{\rm pr}_1^*{\cal F},{\rm pr}_2^!{\cal F})$,
the map (\ref{eqeva}) is the dual of
the map
$\Lambda_X\to \delta^!{\cal H}^*$
corresponding to
$1=u(X,1)$.

\begin{pr}\label{prKS}
Let $c\colon C\to X\times X$ be a closed immersion
and $u\colon c_!\Lambda_C
\to {\cal H}=R{\cal H}om(
{\rm pr}_2^*{\cal F},{\rm pr}_1^!{\cal F})$
be a cohomological correspondence.
Let $i\colon C\cap X\to X$ denote the immersion.
Then 
the characteristic class
$C({\cal F},C,u)\in H^0_{C\cap X}(X,{\cal K}_X)$
is equal to the cohomology class
of the composition 
\begin{equation}
\begin{CD}
i_!\Lambda_{C\cap X}
=\delta^*c_!\Lambda_C
@>{\delta^*(u)}>> \delta^*{\cal H}
@>{\rm(\ref{eqeva})}>>{\cal K}_X.
\end{CD}
\label{eqdel}
\end{equation}
\end{pr}

{\it Proof.}
Let $1\colon \Lambda_X\to \delta^!{\cal H}^*$
be the cohomological correspondence
defined by the identity of ${\cal F}$.
Then the evaluation map
(\ref{eqeva}) is equal to the composition of
$$\begin{CD}
\delta^*{\cal H}
@>{{\rm id}\otimes 1}>>
\delta^*{\cal H}
\otimes \delta^!{\cal H}^*
@>{\rm(\ref{FoxG})}>>
\delta^!({\cal H}
\otimes {\cal H}^*)
@>>> 
\delta^!{\cal K}_{X\times X}={\cal K}_X
\end{CD}$$
where the last arrow is induced
by the pairing {\rm(\ref{eqVerdi})}.
Thus the assertion follows from Lemma \ref{lmcomp}.
\qed

We define a refinement 
of the characteristic class for zero-extensions.
Let $X$ be a scheme over $k$ and
$j_U\colon U\to X$ be an open immersion over $k$.
Let $\delta\colon X
\to X\times X$ and
$\delta_U\colon U
\to U\times U$ be the diagonal maps.
Let ${\cal F}$ be an object of $D_{\rm ctf}(U)$.
We put ${\cal H}=
R{\cal H}om({\rm pr}_2^*{\cal F},
{\rm pr}_1^!{\cal F})$
on $U\times U$
and $\overline{\cal H}=
R{\cal H}om({\rm pr}_2^*j_{U!}{\cal F},
{\rm pr}_1^!j_{U!}{\cal F})$
on $X\times X$.
Since $\overline{\cal H}=
(j_U\times 1)_!
(1\times j_U)_*{\cal H}$,
we have a canonical isomorphism
$j_{U!}\delta_U^*{\cal H}\to
\delta^*\overline{\cal H}$.
Thus the evaluation map
$e\colon \delta_U^*{\cal H}\to {\cal K}_U$
(\ref{eqeva}) on $U$
induces its zero-extension
$j_{U!}e\colon \delta^*\overline{\cal H}\to 
j_{U!}{\cal K}_U$.

\begin{df}\label{dfccc}
Let $X$ be a scheme over $k$ and
$j_U\colon U\to X$ be an open immersion over $k$.
Let $\overline C$ be a closed subscheme of $X\times X$
and we consider the cartesian diagram
$$\begin{CD}
C=&\overline C\cap (U\times U)
@>{j}>> \overline C\\
&@V{c}VV @VV{\bar c}V\\
&U\times U@>>> X\times X.
\end{CD}$$
We assume 
$C=\overline C\cap (X\times U)$.
Let $\bar i\colon {\overline C}\cap X\to X$
denotes the closed immersion.
Let ${\cal F}$ be an object of $D_{\rm ctf}(U)$.
Let 
$u\colon c_!\Lambda \to {\cal H}$
be a cohomological correspondence
and let
$j_!u\colon \bar c_!\Lambda \to \overline {\cal H}$
be the zero-extension.
Then, we define the refinement 
$C_!(j_{U!}{\cal F},\overline C,j_!u)
\in H^0_{\overline C\cap X}(X,j_{U!}{\cal K}_U)$
of the characteristic class
to be the class of the composition
$$\begin{CD}
\bar i_!\Lambda_{\overline C\cap X}
=\delta^*\bar c_!\Lambda_C
@>{\delta^*(j_!u)}>>
\delta^*\overline {\cal H}
@>{j_{U!}e}>>
j_{U!}{\cal K}_U.
\end{CD}$$

If $\overline C$ is the diagonal
and $u\colon {\cal F}\to {\cal F}$
is an endomorphism,
we drop $\overline C$ from the notation and we
write $C_!(j_{U!}{\cal F},j_!u)
\in H^0(X,j_{U!}{\cal K}_U)$.
Further, if $u$ is the identity,
we simply write 
$C_!(j_{U!}{\cal F})$
and call it the refined
characteristic class of $j_{U!}{\cal F}$.
\end{df}

It is clear from Proposition \ref{prKS}
that 
the characteristic class
$C(j_{U!}{\cal F},\overline C,j_!u)
\in H^0_{\overline C\cap X}(X,{\cal K}_X)$
is the image of the refinement
$C_!(j_{U!}{\cal F},\overline C,j_!u)
\in H^0_{\overline C\cap X}(X,j_{U!}{\cal K}_U)$
by the canonical map
$H^0_{\overline C\cap X}(X,j_{U!}{\cal K}_U)
\to H^0_{\overline C\cap X}(X,{\cal K}_X)$.

We give a compatibility with
the pull-back.
Let 
$$\begin{CD}
V@>{j_V}>> Y\\
@VfVV @VV{\bar f}V\\
U@>{j_U}>>X
\end{CD}$$ 
be a cartesian square of schemes
over $k$ 
where the horizontal arrows are
open immersions.
We assume that $U$ and $V$
are smooth over $k$
and $\dim U=\dim V$.
By this assumption,
the map (\ref{f^*f^!})
induces a canonical map
\begin{equation}
\bar f^*j_{U!}{\cal K}_U
=j_{V!}f^*{\cal K}_U
\to j_{V!}f^!{\cal K}_U=
j_{V!}{\cal K}_V
\label{eqf^*j}
\end{equation}
Let $c=(c_1,c_2)\colon 
C\to U\times U$
be a closed immersion.
Let $\overline C$ be the closure of $C$
in $X\times X$
and put 
$\overline C'=(\bar f\times \bar f)^{-1}\overline C$.
Then, the map (\ref{eqf^*j})
induces a map
$\bar f^*\colon H^0_{\overline C\cap X}(X,j_{U!}{\cal K}_U)
\to H^0_{\overline C'\cap Y}(Y,j_{V!}{\cal K}_V)$.

\begin{pr}\label{prpull!}
Let the notation be as above and
assume that 
$C=\overline C\cap (X\times U)$.
Let ${\cal F}$ be an object of $D_{\rm ctf}(U)$
and let 
$u\colon c_2^*{\cal F}\to c_1^!{\cal F}$
be a cohomological correspondence on $C$.
Let $j\colon C\to \overline C$ 
and $j'\colon C'=(f\times f)^{-1}C
\to \overline C'$
denote the open immersions.
Then, we have
$C'=\overline C'\cap (Y\times V)$ and
$$
C_!(j_{V!}f^*{\cal F},\overline C',j'_!(f\times f)^*u)
=\bar f^*C_!(j_{U!}{\cal F},\overline C,j_!u)
$$
in $H^0_{\overline C'\cap Y}(Y,j_{V!}{\cal K}_V)$.
\end{pr}

{\it Proof.}
The equality
$C'=\overline C'\cap (Y\times V)$ is
clear from
$C=\overline C\cap (X\times U)$.
We put
${\cal H}
=R{\cal H}om({\rm pr}_2^*{\cal F},
{\rm pr}_1^!{\cal F})$
and $\overline{\cal H}
=R{\cal H}om({\rm pr}_2^*j_{U!}{\cal F},
{\rm pr}_1^!j_{U!}{\cal F})$
on $U\times U$
and on $X\times X$
respectively.
We also put
${\cal H}'
=R{\cal H}om({\rm pr}_2^*f^*{\cal F},
{\rm pr}_1^!f^*{\cal F})$
on $V\times V$ and
$\overline{\cal H}'
=R{\cal H}om({\rm pr}_2^*j_{V!}f^*{\cal F},
{\rm pr}_1^!j_{V!}f^*{\cal F})$
on $Y\times Y$.
We consider the commutative diagram
$$\begin{CD}
j_!u\in&
H^0_{\overline C}(X\times X,
\overline{\cal H})
@>{\delta_X^*}>>
H^0_{\overline C\cap X}(X,j_{U!}\delta_U^*{\cal H})
@>{j_{U!}e}>>
H^0_{\overline C\cap X}(X,j_{U!}{\cal K}_U)\\
&@V{(\bar f\times \bar f)^*}VV 
@VV{\bar f^*}V
@VV{\bar f^*}V\\
j'_!(f\times f)^*u
\in&
H^0_{\overline C'}(Y\times Y,
\overline{\cal H}')
@>{\delta_Y^*}>>
H^0_{\overline C'\cap Y}(Y,j_{V!}\delta_V^*{\cal H}')
@>{j_{V!}e}>>
H^0_{\overline C'\cap Y}(Y,j_{V!}{\cal K}_V).
\end{CD}$$
By the uniqueness Lemma \ref{lmext},
we have
$(\bar f\times \bar f)^*j_!u
=j'_!(f\times f)^*u$.
Thus the assertion follows.
\qed

\begin{cor}\label{corpull}
If $X=U$ and $Y=V$, we have
$$C(f^*{\cal F},C',(f\times f)^*u)=f^*C({\cal F},C,u).$$
\end{cor}

We apply 
Proposition \ref{prpull!}
to a finite \'etale Galois covering.
We keep the notation in
Proposition \ref{prpull!}
and assume that 
$V\to U$
is a finite \'etale Galois covering
of Galois group $G$.
We assume further that
$C=U$ is the diagonal
and hence $u\colon {\cal F}\to {\cal F}$
is an endomorphism.
For $\sigma\in G$,
let $\Gamma_\sigma\subset V\times V$
be the graph of $\sigma\colon V\to V$
and let $\sigma^*\colon \sigma^*f^*{\cal F}\to 
(f\circ \sigma)^*{\cal F}=f^*{\cal F}$
denote the canonical map.
We consider
the composition
$f^*(u)\circ \sigma^*\colon 
\sigma^*f^*{\cal F}\to f^*{\cal F}$
as a cohomological correspondence
on $f^*{\cal F}$ on the graph
$\Gamma_\sigma\subset V\times V$.

\begin{cor}\label{corpull!}
Let the notation be as in Proposition \ref{prpull!}.
We assume that $C=\Delta_U\subset U\times U$
and $u\colon {\cal F}\to {\cal F}$
is an endomorphism of ${\cal F}$.
We assume further that 
$V\to U$
is a finite \'etale Galois covering
of Galois group $G$.
For $\sigma\in G$,
let $j_\sigma\colon \Gamma_\sigma\to
\overline{\Gamma_\sigma}$
be the open immersion to the closure
in $Y\times Y$.
Then, we have
$$\bar f^*C_!(j_{U!}{\cal F},j_!u)
=\sum_{\sigma\in G}
C_!(j_{V!}f^*{\cal F},\overline{\Gamma_\sigma},
j_{\sigma!}(f^*(u)\circ \sigma^*))$$
in $H^0(Y,j_{V!}{\cal K}_V)$.
\end{cor}

{\it Proof.}
By Proposition \ref{prpull!},
it suffices to show
$C_!(j_{V!}f^*{\cal F},\overline C',j'_!(f\times f)^*u)
=\sum_{\sigma\in G}
C_!(j_{V!}f^*{\cal F},\overline{\Gamma_\sigma},
j_{\sigma!}(f^*(u)\circ \sigma^*))$.
Since $V\to U$
is an \'etale Galois covering,
we have
$(f\times f)^*u
=\sum_{\sigma\in G}
f^*(u)\circ \sigma^*$.
Since
$\overline{\Gamma_\sigma}
\subset \overline C'$,
we have 
$C_!(j_{V!}f^*{\cal F},\overline C',
j'_!(f^*(u)\circ \sigma^*))
=C_!(j_{V!}f^*{\cal F},\overline{\Gamma_\sigma},
j_{\sigma!}(f^*(u)\circ \sigma^*))$
in $H^0(Y,j_{V!}{\cal K}_V)$ and
the assertion follows.
\qed

\begin{pr}\label{prpink}
We consider a commutative diagram
\begin{equation}
\xymatrix{
& C\ar[dl]\ar[dddl]^{j}\ar[rr]^{\!\!\!\! c}
&&U\times U\ar[dl]\ar[dddl]^{j_U\times j_U}&&
U\ar[dl]\ar[dddl]^{j_U}\ar[ll]_{\quad\delta}\\
{\overline C}'\ar[dd]_{g}\ar[rr]^{\!\!\!\!\!\! \bar c'}&&
(X\times X)'
\ar[dd]_{f}&&
X\ar[dd]_{\|}\ar[ll]_{\qquad \qquad \bar \delta'}
&\\\\
{\overline C}\ar[rr]^{\!\!\!\! \bar c}&&
X\times X&&
X\ar[ll]_{\quad\bar \delta}&}
\label{eqpink2}
\end{equation}
of schemes over $k$.
We assume 
that the horizontal arrows are closed immersions, 
that the six slant arrows are open immersions,
that the vertial arrows $f$ and $g$
are proper and that 
$f$ is an isomorphism on 
$(X\times U)\cup (U\times X)$.
We also assume that
the four slant parallelograms are cartesian
and 
$C=
\overline C\times_{X\times X}(X\times U)$.

Let ${\cal F}$ be an object of
$D_{\rm ctf}(U)$.
We put 
${\cal H}=R{\cal H}om(
{\rm pr}_2^*{\cal F},{\rm pr}_1^!{\cal F})$
and 
${\cal H}^*=R{\cal H}om(
{\rm pr}_1^*{\cal F},{\rm pr}_2^!{\cal F})$
on $U\times U$.
We also put
$\overline {\cal H}=R{\cal H}om(
{\rm pr}_2^*j_!{\cal F},{\rm pr}_1^!j_!{\cal F})$
and 
$\overline {\cal H}^*=R{\cal H}om(
{\rm pr}_1^*j_!{\cal F},{\rm pr}_2^!j_!{\cal H})$
on $X\times X$.

Let
$\overline {\cal H}'$
and 
$\overline {\cal H}^{*\prime}$
be objects of
$D_{\rm ctf}((X\times X)')$
and $a\colon 
\overline {\cal H}|_{U\times Y}\to
\overline {\cal H}'|_{U\times Y}$
and
$b\colon 
\overline {\cal H}^*|_{X\times V}\to
\overline {\cal H}^{*\prime}|_{X\times V}$
be isomorphisms such
that the unique maps
$\bar a\colon \overline {\cal H}
\to f_*\overline {\cal H}'$
and $\bar b\colon \overline {\cal H}^*
\to f_*\overline {\cal H}^{*\prime}$
inducing $a$ and $b$
in Proposition {\rm \ref{lmpink}.1}
are isomorphisms.
Let $\overline {\cal H}'
\otimes 
\overline {\cal H}^{*\prime}
\to {\cal K}_{(X\times X)'}$
be a map making the diagram
{\rm (\ref{eqpinkp})} commutative.

Let
$u\colon c_!\Lambda\to {\cal H}$
be a cohomological correspondence.
Let $\bar u'
\colon \bar c'_!\Lambda\to \overline {\cal H}'$
and $\bar 1'
\colon \bar \delta'_!\Lambda\to \overline {\cal H}^{\prime *}$
be maps making
the diagrams {\rm (\ref{eqlmpink1})}
commutative for $\bar d=\bar \delta\colon X\to X\times X$
and $v=1$.

1. We have
$$C(j_!{\cal F},\overline C,j_!u)=
\langle \bar u',\bar 1'\rangle$$
in $H^0(\overline C\cap X,
{\cal K}_{\overline C\cap X})
=H^0_{\overline C\cap X}(X,
{\cal K}_X)$.

2. Further, we regard
$\bar u'\in
H^0_{\overline C'}(
(X\times X)',\overline{\cal H}')$
and let $\bar \delta^{\prime*}\colon 
H^0_{\overline C'}(
(X\times X)',\overline{\cal H}')
\to
H^0_{\overline C'\cap X}(
X,\bar \delta^{\prime*}\overline{\cal H}')$
be the pull-back.
Let $\widetilde{\rm Tr}\colon 
\bar \delta^{\prime*}\overline{\cal H}'
\to K_X$
be the pairing
$\bar \delta^{\prime*}\overline{\cal H}'
\otimes
\bar \delta^{\prime!}
\overline{\cal H}^{*\prime}
\to 
\bar \delta^{\prime!}
K_{(X\times X)'}=K_X$
{\rm(\ref{FoxG})}
induced by
$\overline{\cal H}'
\otimes
\overline{\cal H}^{*\prime}
\to 
K_{(X\times X)'}$
evaluated at $\bar 1'$.
Then, we have
$$C(j_!{\cal F},\overline C,j_!u)=
\widetilde{\rm Tr}\
\bar \delta^{\prime*} \bar u'$$
in $H^0(\overline C\cap X,
{\cal K}_{\overline C\cap X})
=H^0_{\overline C\cap X}(X,
{\cal K}_X)$.
\end{pr}

{\it Proof.}
1. By Corollary \ref{corext},
we have $1=j_!1$.
Hence, we have
$C(j_!{\cal F},\overline C,j_!u)=
\langle j_!u,j_!1\rangle$.
Since the assumptions in
Proposition \ref{lmpink} are satisfied,
we have
$\langle j_!u,j_!1\rangle=
f_*\langle \bar u',\bar 1'\rangle$.
Since the restriction of $f$
to $X$ is the identity,
the assertion follows.

2. We have
$\langle \bar u',\bar 1'\rangle=
\widetilde{\rm Tr}\
\delta^{\prime*} \bar u'$
by Lemma \ref{lmcomp}.
Thus, it follows from 1.
\qed

\subsection{Characteristic class 
and log blow-up}

In this subsection,
$X$ denotes a smooth scheme  
of dimension $d$ over $k$
and $D=\bigcup_{i=1}^mD_i\subset X$
denotes a divisor with simple normal crossings.
Let $U\subset X$ be the complement of $D$
and $j\colon U\to X$ be the open immersion.

Let $(X\times X)'\to X\times X$ be the blow-up
at $D_i\times D_i$ for
$i=1,\ldots,m$.
It is defined by the product of
the ideal sheaves
$I_{D_i\times D_i}\subset O_{X\times X}$.
Let $p_1,p_2\colon (X\times X)'\to X$
be the projections.
For $i=1,\ldots,m$,
let $(D_i\times X)'$ be the proper transform of
$D_i\times X$ and $(X\times D_i)'$ 
be that of $X\times D_i$. 
We put
$D^{(1)\prime}=
\bigcup_i
(D_i\times X)'$ and
$D^{(2)\prime}
=\bigcup_i
(X\times D_i)'$.
We define open subschemes
$(U\times X)'$ and $(X\times U)'$
of $(X\times X)'$
to be the complements of
$D^{(1)\prime}$ and of
$D^{(2)\prime}$ 
respectively.

The intersection
$(X\times X)^\sim
=(U\times X)'\cap (X\times U)'$
is the maximum open subscheme
of $(X\times X)'$
where
$p_1^*D_i=p_2^*D_i$ 
for each $i=1,\ldots,m$.
We obtain a commutative diagram
\begin{equation}
\xymatrix{
 &(X\times U)'\ar[dl]_{j'_2}&&(X\times 
X)^\sim \ar[dl]_{k'_2}\ar[ll]_{k'_1}\\
(X\times X)'\ar[dd]_f&&(U\times X)' \ar[ll]_{\qquad \qquad j'_1}&\\
 &X\times U\ar[dl]_{j_2}\ar[uu]&&
U\times U\ar[dl]_{k_2}\ar[ll]_{\qquad \qquad k_1}\ar[uu]_{\tilde j}\\
X\times X&&U\times X.\ar[ll]_{j_1}\ar[uu]&}
\end{equation}
All the arrows
except the blow-up
$f\colon (X\times X)'\to X\times X$
are open immersions.
The four faces consisting of open immersions
are cartesian.
The diagonal map $\delta_X=\delta\colon X=\Delta_X\to X\times X$
induces the 
log diagonal map $\delta'_X=\delta'\colon X
=\Delta_X^{\log}\to 
(X\times X)^\sim\subset
(X\times X)'$.

Let ${\cal F}$ be a locally constant
sheaf of free $\Lambda$-modules
of finite rank on $U$
tamely ramified along
$D=X\setminus U$.
We put ${\cal H}_0={\cal H}om({\rm pr}_2^*{\cal F},
{\rm pr}_1^*{\cal F})$
and ${\cal H}=R{\cal H}om({\rm pr}_2^*{\cal F},
{\rm pr}_1^!{\cal F})
={\cal H}_0(d)[2d]$ 
on $U\times U$.
We put $\overline {\cal H}=
R{\cal H}om({\rm pr}_2^*j_!{\cal F},
{\rm pr}_1^!j_!{\cal F})=
j_{1!}Rk_{2*}{\cal H}=
Rj_{2*}k_{1!}{\cal H}$ on $X\times X$.
We put $\widetilde {\cal H}_0=\tilde j_*{\cal H}_0$
on $(X\times X)^\sim$.
Here and in the rest of this subsection,
the symbol $j_*$ for an open immersion $j$
denotes the usual direct image
and is not an abbreviation of
$Rj_*$.
We also put $\widetilde {\cal H}=\widetilde {\cal H}_0(d)[2d]$.
The trace map
${\rm Tr}\colon \delta_U^*{\cal H}_0={\cal E}nd({\cal F})
\to \Lambda_U$
is extended uniquely to a map
\begin{equation}
\begin{CD}
\widetilde {\rm Tr}\colon \delta^{\prime *}\widetilde {\cal H}_0
=j_*{\cal E}nd({\cal F})@>>>
\Lambda_X.
\end{CD}
\label{eqtamej}
\end{equation}
We define $\overline {\cal H}'=
j'_{1!}Rk'_{2*}\widetilde {\cal H}=
Rj'_{2*}k'_{1!}\widetilde {\cal H}$
on $(X\times X)'$.

\begin{df}\label{dftame}
Let $j\colon U\to X$
and ${\cal F}$ on $U$ be as above.
Let $C$
be a closed subscheme of $U\times U$
and let $\overline C'$ and $\widetilde C
=\overline C'\cap (X\times X)^\sim$
be the closures of $C$
in $(X\times X)'$ and
in $(X\times X)^\sim$
respectively.

1. 
We say that a closed subscheme
$C\subset U\times U$ is non-expanding
with respect to $X$
if we have
$\widetilde C=\overline C'\cap (X\times U)'$.

2. 
Let  
$C\subset U\times U$ be a closed subscheme
that is non-expanding with respect to $X$.
We call the image 
of the injection
$\Gamma(\widetilde C,
\widetilde {\cal H}_0|_{\widetilde C})
\to 
\Gamma(C,{\cal H}_0|_C)$
the tame part of
$\Gamma(C,{\cal H}_0|_C)
=Hom(c_2^*{\cal F},c_1^*{\cal F})$.

3. Let
$\gamma\colon c_2^*{\cal F}\to c_1^*{\cal F}$
be a map in the tame part.
For a cycle class $\widetilde \Gamma\in CH_d(\widetilde C)$,
let $\tilde u=\tilde u(\widetilde \Gamma,\gamma)
\in H^0_{\widetilde C}((X\times X)^\sim,\widetilde {\cal H})$
denote the cup-product of
${\rm cl}(\widetilde \Gamma)\in 
H^0_{\widetilde C}((X\times X)^\sim,\Lambda(d)[2d])$
with
$\gamma\in {\Gamma}(\widetilde C,\widetilde {\cal H}_0)$.
Let $\widetilde {\rm Tr}(\delta^{\prime *}\gamma)
\in \Gamma(\widetilde C\cap X,\Lambda)$
denote the image of
$\delta^{\prime *}\gamma\in 
\Gamma(\widetilde C\cap X,\widetilde {\cal H}_0)$
by the map induced by $\widetilde {\rm Tr}$
{\rm (\ref{eqtamej})}.
\end{df}

If a closed subscheme $C\subset U\times U$
is non-expanding with respect to $X$,
the closure $\overline C
\subset X\times X$
of $C$
satisfies
$\overline C\cap (X\times U)=
\widetilde C\cap (X\times U)=C$.
For $\widetilde \Gamma\in CH_d(\widetilde C)$
and $\gamma\colon c_2^*{\cal F}\to c_1^*{\cal F}$
in the tame part,
the element $\tilde u=\tilde u(\widetilde \Gamma,\gamma)
\in H^0_{\widetilde C}((X\times X)^\sim,\widetilde {\cal H})$
is an extension of
$u=u(\Gamma,\gamma)
\in H^0_C(U\times U,{\cal H})$
defined by (\ref{equgam}).

The following lemma gives
an example of a non-expanding closed subscheme.

\begin{lm}\label{lmtameC}
{\rm (\cite{KS} Proposition 1.1.6)}
Let $j\colon U\to X$ be as above.
We consider a cartesian diagram
\begin{equation}
\begin{CD}
U@>>> X\\
@VVV @VVfV\\
V@>>> Y
\end{CD}
\label{eqcart}
\end{equation}
of morphisms of schemes
over $k$.
We assume $V$ is the complement of
a Cartier divisor $B$ of $Y$.
Then,
the closed subscheme $C=U\times_VU
\subset U\times U$ is non-expanding
with respect to $X$.
\end{lm}

In the tame case,
the characteristic class is
computed as follows.

\begin{pr}\label{prtame}
Let the notation be as above.
Let ${\cal F}$ be a smooth
$\Lambda$-sheaf on $U=X\setminus D$ tamely ramified
along $D$.
Let $C\subset U\times U$ be a
closed subscheme that is non-expanding with
respect to $X$
and 
$\overline C\subset X\times X$,
$\overline C'\subset (X\times X)'$ and
$\widetilde C\subset (X\times X)^\sim$
be the closures.

1. 
The restriction map
\begin{equation}
H^0_{\overline C'}
((X\times X)',
\overline{\cal H}')
\to
H^0_{\widetilde C}
((X\times X)^\sim,
\widetilde{\cal H})
\label{eqtameex}
\end{equation}
is an isomorphism.

2. 
Let $\gamma\colon c_2^*{\cal F}
\to c_1^*{\cal F}$ be
an element in the tame part
$\Gamma(\widetilde C,
\widetilde {\cal H}_0|_{\widetilde C})
\subset 
\Gamma(C,{\cal H}_0|_C)
=Hom(c_2^*{\cal F}, c_1^*{\cal F})$.
Let $\widetilde \Gamma \in CH_0(\widetilde C)$
be an algebraic cycle class 
and $\Gamma\in CH_0(C)$
be the restriction.
Then, we have
\begin{equation}
C(j_!{\cal F},\overline C,j_!u(\Gamma,\gamma))=
{\rm cl}(\widetilde \Gamma,\Delta_X^{\log})_{
(X\times X)^\sim}\cdot 
\widetilde {\rm Tr}\ \delta^{\prime *}\gamma
\label{eqtame}
\end{equation}
in $H^{2d}_{\overline C\cap X}(X,\Lambda(d))=
H^0_{\overline C\cap X}(X,{\cal K}_X)$.
The right hand side
is the product of the class
${\rm cl}(\widetilde \Gamma,\Delta_X^{\log})_{
(X\times X)^\sim}
\in 
H^0_{\overline C\cap X}(X,{\cal K}_X)$
of the intersection product and
$\widetilde {\rm Tr}\ \delta^{\prime *}\gamma
\in \Gamma(\overline C\cap X,\Lambda)$
defined in Definition {\rm \ref{dftame}.3}.
\end{pr}

{\it Proof.}
1. 
Since $\overline C'\cap 
(X\times U)'=\widetilde C$,
we have an isomorphism
$$\begin{array}{cl}
&H^0_{\overline C'}
((X\times X)',
\overline{\cal H}')
=
H^0_{\overline C'}
((X\times X)',
Rj'_{2*}k'_{1!}\widetilde {\cal H})\\
&\\
\to&
H^0_{\overline C'\cap 
(X\times U)'}
((X\times U)',
k'_{1!}\widetilde {\cal H})
=
H^0_{\widetilde C}
((X\times X)^\sim,
\widetilde {\cal H})
\end{array}
$$
as in (\ref{eqzeroex}).

2. Similarly as 
$\overline{\cal H}'
=j'_{1!}Rk'_{2*}\tilde j_*{\cal H}_0(d)[2d]$,
we put
${\cal H}_0^*=
{\cal H}om({\rm pr}_1^*{\cal F},
{\rm pr}_2^*{\cal F})$
on $U\times U$,
$\widetilde {\cal H}^*=
\tilde j_*{\cal H}_0^*(d)[2d]$
on $(X\times X)^\sim$,
and
$\overline{\cal H}^{\prime*}
=j'_{2!}Rk'_{1*}\widetilde {\cal H}^*
=Rj'_{1*}k'_{2!}\widetilde {\cal H}^*$
on $(X\times X)'$.
By 1,
we regard
$\tilde u(\tilde \Gamma,\gamma)
\in H^0_{\widetilde C}((X\times X)^\sim,
\widetilde {\cal H})$
and
$\tilde u(X,1)
\in H^0_X((X\times X)^\sim,
\widetilde {\cal H}^{*\prime})$
as elements of
$H^0_{\overline C'}((X\times X)',
\overline{\cal H}')$
and
$H^0_X((X\times X)',
\overline{\cal H}^{*\prime})$
respectively.
We consider a pairing
$\overline{\cal H}'
\otimes
\overline{\cal H}^{\prime*}
=(j'_1k'_2)_!
(\tilde j_*{\cal H}_0
\otimes
\tilde j_*{\cal H}_0^*)
(2d)[4d]
\to {\cal K}_{(X\times X)'}$
induced by the pairing
${\cal H}_0
\otimes
{\cal H}_0^*
\to \Lambda_{U\times U}\colon 
a\otimes b\mapsto {\rm Tr}(b\circ a)$.

By Lemma \ref{lmFP} below,
the assumptions 
in Proposition \ref{prpink}
are satisfied.
We claim that the map
$\widetilde{\rm Tr}\colon 
\bar \delta^{\prime*}\overline{\cal H}'
=j_*{\cal E}nd({\cal F})(d)[2d]
\to K_X=\Lambda_X(d)[2d]$ 
in Proposition \ref{prpink}.2
is induced by the map
$\widetilde{\rm Tr}$ defined in
(\ref{eqtamej}).
On $U$, it is shown in the proof
of Proposition \ref{prKS}.
Since the map
$\widetilde{\rm Tr}\colon 
\delta^{\prime*}\widetilde {\cal H}_0
=j_*{\cal E}nd({\cal F})
\to \Lambda_X$ 
is the unique map extending
${\rm Tr}\colon 
\delta^* {\cal H}_0
={\cal E}nd({\cal F})
\to \Lambda_U$,
it induces the map 
$\widetilde{\rm Tr}\colon 
\bar \delta^{\prime*}\overline{\cal H}'
\to K_X$.
Thus, 
applying 
Proposition \ref{prpink}.2,
we obtain
$C(j_!{\cal F},\overline C,j_!u)
=
\widetilde{\rm Tr}\
\delta^{\prime*} \tilde u(\tilde \Gamma,\gamma)
=
\delta^{\prime*} ({\rm cl}(\tilde \Gamma))\cdot
\widetilde{\rm Tr}\
\delta^{\prime*} \gamma$.
By the compatibility (\ref{eqcycpull}),
we have
$\delta^{\prime*} ({\rm cl}(\tilde \Gamma))=
{\rm cl}(\tilde \Gamma,\Delta_X^{\log})
_{(X\times X)^\sim}$.
Thus the assertion follows.
\qed

\begin{lm}\label{lmFP}
{\rm (cf. \cite{fal},
\cite{Pink}, \cite{KS} Lemma 2.2.1)}
The canonical maps
$\overline {\cal H}
\to f_*\overline {\cal H}'$
and $\overline {\cal H}^*
\to f_*\overline {\cal H}^{*\prime}$
are isomorphisms.
\end{lm}
{\it Proof.}
It suffices to show that
the stalk of
$f_*\overline {\cal H}'$
is 0 at an arbitrary
geometric point of $\bigcup_iD_i\times D_i$.
Let $x\in \bigcup_iD_i\times D_i$
and $D_1,\ldots,D_m$
be the components of $D$.
Then, since the question is \'etale local
and since it is reduced to
the case where $\Lambda$ is a field,
we may assume ${\cal F}$
is the tensor product
${\cal F}_1\otimes \cdots\otimes {\cal F}_m$
where ${\cal F}_i$
is the extension by zero
of a smooth sheaf of rank 1
on the complement $X\setminus D_i$
for each $i=1,\ldots,r$.

Let $(X\times X)'_i$
be the blow-up of $X\times X$
at $D_i\times D_i$.
Then, $(X\times X)'$
is the fiber product
of $(X\times X)'_i$
over $X\times X$.
Hence, by the K\"unneth formula,
it is reduced to the case $r=1$ and $D=D_1$.

The fiber of
$(X\times X)^\sim\subset (X\times X)'\to
X\times X$ at a point in $D_1\times D_1$
is isomorphic to ${\mathbb G}_m
\subset {\mathbf P}^1$
and the restriction of 
${\cal H}_0$ to ${\mathbb G}_m$
is smooth and is tamely ramified
along $0$ and $\infty$.
Hence, by the proper base change theorem,
it is reduced to showing
that $H^q({\mathbb A}^1_{\overline k},j_!{\cal G})=0$
for $q=0,1,2$
for a smooth sheaf ${\cal G}$
on ${\mathbb G}_m$ tamely ramified at $0$ and $\infty$
where $j\colon {\mathbb G}_m\to
{\mathbb A}^1$
is the open immersion.

It is clear for $q=0$
and follows by duality for $q=2$.
The smooth sheaf ${\cal G}$ is trivialized by
the covering $[m]\colon {\mathbb G}_m
\to {\mathbb G}_m$
for an integer $m$ invertible in $k$.
Thus, the case $q=1$ is reduced to the case
where ${\cal G}$ is constant.
Hence it follows from the acyclicity of
${\mathbb A}^1$.
\qed

We put $C_U=(-1)^dc_d(\Omega^1_{X/k}(\log D))
\in H^{2d}(X,{\mathbb Z}_\ell(d))$.
Since the conormal sheaf
$N_{X/(X\times X)^\sim}$
is canonically isomorphic to
$\Omega^1_{X/k}(\log D)$,
we have $C_U=(\Delta_X^{\log},
\Delta_X^{\log})_{(X\times X)^\sim}$.

\begin{cor}\label{cortame}
1. Assume that 
${\cal F}$ is tamely ramified
along $D$
and that $C=U$ is the diagonal.
Then every endomorphism 
of ${\cal F}$ is in the tame part and
we have
${\rm End}_U({\cal F})^\sim
={\rm End}_U({\cal F})$.
For an endomorphism $\gamma\in 
{\rm End}_U({\cal F})$, 
we have
$$C(j_!{\cal F},j_!\gamma)=
{\rm Tr}\ \gamma\cdot C_U$$
in $H^{2d}(X,\Lambda(d))$.
In particular, for ${\cal F}=\Lambda$
and $\gamma={\rm id}$,
we have
$$C(j_!\Lambda)=
C_U.$$

2. Assume $C$ is non-expanding
and is of dimension $d$.
Then, we have
\begin{equation}
C(j_!\Lambda,C,u(C,1))=
(\widetilde C,\Delta_X^{\log})_{(X\times X)^\sim}.
\label{eqcort}
\end{equation}
\end{cor}

If $X$ is proper,
taking the trace of
(\ref{eqcort}),
we recover the Lefschetz trace formula
$${\rm Tr}(C^*,H^*_c(U_{\overline k},\Lambda))=
{\rm Tr}(\widetilde C,\Delta_X^{\log})_{(X\times X)^\sim}$$
for an open variety
\cite{KS} Theorem 2.3.4. In particular,
if $C=U$, we have
$$\chi_c(U)=
{\rm deg}\ C_U=
{\rm deg}\ (\Delta_X^{\log},
\Delta_X^{\log})_{(X\times X)^\sim}.$$

\section{Characteristic class and Swan class}

In this section,
$\ell$ denotes a prime number
invertible on a noetherian scheme $X$
and 
$E$ denotes
a finite extension of
${\mathbb Q}_\ell$.
Let $O$ be the integer ring of $E$
and
$F=O/\lambda$ be
the residue field of $E$.
For a constructible $E$-sheaf ${\cal F}$ on $X$,
${\cal F}_O$ denotes an $O$-lattice
and 
${\cal F}_n$
denotes the reduction
${\cal F}_O\otimes_OO/\lambda^n$.
We put $\overline {\cal F}=
{\cal F}_1$.

For a constructible $E$-sheaf ${\cal F}$
on a scheme $X$ over $k$
and a cohomological correspondence
$u$ on $C\subset X\times X$,
the characteristic class
$C({\cal F},C,u)\in H^0(X,{\cal K}_X)$
is defined as
$C({\cal F}_n,C,u_n)_n
\in \varprojlim_n H^0(X,Rf^!O/\lambda^n)\otimes_OE$.
It inherits the properties
established in the previous
section for torsion coefficients.
For a smooth $E$-sheaf ${\cal F}$
on a smooth dense open subscheme
$U\subset X$,
we will prove that the difference
$C(j_!{\cal F})-
{\rm rank}\ {\cal F}\cdot
C(j_!E)$
is given by the Swan class ${\rm Sw}({\cal F})$
defined in \cite{KS},
under a certain assumption
formulated using a notion introduced in \S3.1.
The definition of the Swan class ${\rm Sw}({\cal F})$
will be recalled in \S3.2.

\subsection{$\ell$-adic sheaves potentially of Kummer type}

\begin{df}\label{dfKum}
Let $X$ be a scheme and $U\subset X$
be an open subscheme.

1. Let $D_i\subset X,i\in I$ be a finite family of 
Cartier divisors
such that $U\cap D_i=\emptyset$ for each $i\in I$.
Then, we say a locally constant constructible sheaf 
${\cal F}$ of sets on $U$ is 
of Kummer type with respect to $D_i,i\in I$
if the following condition is satisfied.
\begin{enumerate}\item[{\rm (Kum)}]
For each $x\in X$,
there exist an \'etale neighborhood
$W\to X$ of $x$,
bases $t_i$ of the invertible ideals $O(-D_i)|_W$
and an integer $m$ invertible on $W$
such that the pull-back of ${\cal F}$
is constant on
$U\times_XW[T_i,i\in I]/(T_i^m-t_i,i\in I)$.
\end{enumerate}

We say that a finite \'etale scheme $V$ over $U$
is of Kummer type with respect to $D_i,i\in I$
if the locally constant sheaf on $U$ represented by $V$
is of Kummer type with respect to $D_i,i\in I$.

2. We say that a locally constant constructible sheaf 
${\cal F}$ of sets
on $U$
is of Kummer type with respect to $X$
if there exists a finite family of 
Cartier divisors
$D_i\subset X,i\in I$
such that ${\cal F}$ is 
of Kummer type with respect to $D_i,i\in I$.
We say that a smooth $E$-sheaf ${\cal F}$
on $U$ is of Kummer type
with respect to $X$
if there exists a finite family of 
Cartier divisors
$D_i\subset X,i\in I$
such that ${\cal F}_n$ is 
of Kummer type with respect to $D_i,i\in I$
for every $n\ge 0$.

3.
We say that a smooth $E$-sheaf ${\cal F}$
on $U$ is potentially of Kummer type
with respect to $X$
if there exists a cartesian diagram
$$\begin{CD}
V@>>> Y\\
@V{f}VV @VV{\bar f}V\\
U@>>> X
\end{CD}$$
satisfying the following properties:
The map $\bar f\colon Y\to X$ is proper surjective,
$f\colon V\to U$ is finite \'etale
and $f^*{\cal F}$ is of Kummer type with respect to $Y$.
\end{df}

For a geometric point $\bar x$ of
a normal noetherian scheme $X$,
let $X_{\bar x}$ be the strict localization of $X$
at $\bar x$
and let $I_{\bar x}=\pi_1(U\times_XX_{\bar x},\bar \eta_x)$ be the 
inertia group at $\bar x$
where $\bar\eta_{\bar x}$ is a geometric
point of $U_{\bar x}=U\times_XX_{\bar x}$.

\begin{lm}\label{lmKum}
Let $X$ be a normal noetherian scheme
and $U\subset X$ be an open subscheme.
We consider the following conditions on
a locally constant constructible sheaf 
${\cal F}$ of sets on $U$.

(1) ${\cal F}$ is of Kummer type with respect to $X$.

(2) For each geometric point $\bar x$ of $X$,
the image of
$I_{\bar x}\to {\rm Aut}({\cal F}_{\bar\eta_{\bar x}})$
is a finite abelian group of order invertible at $\bar x$.

(3) For each geometric point $\bar x$ of $X$
such that $O_{X,\bar x}$ is a discrete valuation ring,
the image of
$I_{\bar x}\to {\rm Aut}({\cal F}_{\bar\eta_{\bar x}})$
is a finite group of order invertible at $\bar x$.

Then, we have (1) $\Rightarrow$
(2) $\Rightarrow$
(3). 

If $X$ is a regular noetherian scheme
and $U$ is the complement of a divisor
with simple normal crossings,
then the three conditions
(1), (2) and (3) are equivalent. 
\end{lm}

{\it Proof.}
(1) $\Rightarrow$ (2) 
Let ${\cal F}$ be 
a locally constant constructible sheaf 
of sets on $U$
that is of Kummer type with respect to $D_i,i\in I$.
Let the notation be as in the condition
(Kum) and take a geometric point $\bar x$ above $x$.
We define a map
\begin{equation}
I_{\bar x}\to \prod_{i\in I}\mu_m
\label{eqiner}
\end{equation}
by sending $\sigma$ to
$(\sigma(\sqrt[m]{t_i})/\sqrt[m]{t_i})_i$.
Then the image
$I_{\bar x}\to {\rm Aut}({\cal F}_{\bar\eta_{\bar x}})$
is a quotient of the image of
$I_{\bar x}\to \prod_i\mu_m$.

(2) $\Rightarrow$ (3) Clear.

(3) $\Rightarrow$ (1) 
We assume
$X$ is regular
and $D=X\setminus U$
is a divisor with simple normal crossings.
Then
${\cal F}$ is of Kummer type
with respect to
the family
$D_i,i\in I$
of the irreducible components
of $D$
by \cite{SGA1}
Proposition 5.2 (Lemme d'Abhyankhar absolu).
\qed

\begin{cor}\label{corKum}
Let $X$ be a normal scheme over $k$
and $U\subset X$ be a dense open subscheme
smooth over $k$.
If $\dim X\le 2$,
an arbitrary smooth $E$-sheaf ${\cal F}$
on $U$ is potentially of 
Kummer type with respect to $X$.
\end{cor}

{\it Proof.}
Clear from Lemma \ref{lmKum}
and the strong resolution of singularities
for surfaces.
\qed

{\it Remark.}
Let $U\subset X$ be a dense open 
subscheme of a normal scheme $X$
over $k$
and ${\cal F}$ be a locally constant constructible sheaf 
of sets on $U$. We consider the 
following
conditions:
\begin{enumerate}
\item[(1$'$)] There exists a proper modification $X' \to X$ such that 
${\cal F}$ is of Kummer type
with respect to $X'$.
\item[(2$'$)] For each geometric point $\bar x$ of $X$, 
the image of the inertia $I_{\bar x}\to {\rm Aut} {\cal F}_{\bar x}$ 
has order invertible in $k$. 
\end{enumerate}
Kato has told
the authors that (2$'$) implies (1$'$). 
For the other implication,
there is a counterexample.

\medskip

We consider a geometric construction
generalizing the log blow-up
studied in the previous section.
Assume $k$ is perfect and let
$Y$ be a normal scheme over $k$.
Let $V\subset Y$ be an open subscheme
and $j\colon V\to Y$ denote the open immersion.
Let $D_i,i\in I$
be a finite family
of Cartier divisors of $Y$
such that $V\cap D_i=\emptyset$.
Let $(Y\times Y)'$
be the blow-up of $Y\times Y$
at $D_i\times D_i$ for each $i\in I$,
$\overline{(Y\times Y)'}$
be its normalization
and 
let $p_1,p_2\colon \overline{(Y\times Y)'}\to Y$
denote the projections.
Let
$(Y\times Y)^\sim\subset \overline{(Y\times Y)'}$
be the maximum open subscheme 
where $p_1^*D_i=
p_2^*D_i$ for each $i\in I$.
Let $\tilde j\colon 
V\times V\to (Y\times Y)^\sim$
denote the open immersion.

We consider an endomorphism $g\colon Y\to Y$
such that $g^{-1}(V)=V$.
Let $\widetilde Y_g
\subset Y$ be the maximum open subscheme
where $D_i=g^*D_i$ for each $i\in I$.
Then, by the universality of blow-up
and by the assumption that
$Y$ is normal,
the immersion
$\gamma=(1,g)\colon Y\to Y\times Y$
induces an immersion
$\tilde \gamma\colon \widetilde Y_g\to (Y\times Y)^\sim$.
If $g$ is the 
identity of $Y$,
we obtain an immersion
$\tilde \delta\colon Y=\widetilde Y_{\rm id}\to 
(Y\times Y)^\sim$.
Thus, we have a commutative diagram
\begin{equation}
\begin{CD}
V@>{\gamma}>> V\times V @<{\delta}<<V\\
@V{j_g}VV @V{\tilde j}VV @VV{j}V\\
\widetilde Y_g@>{\tilde \gamma}>>
(Y\times Y)^\sim@<{\tilde \delta}<<Y.
\end{CD}
\label{eqjg}
\end{equation}
For a geometric point $\bar y$ of $Y$,
let $Y_{\bar y}$ denote the strict localization.
Taking a geometric point 
$\bar \eta$ of $V_{\bar y}=V\times_YY_{\bar y}$,
we define the inertia group by
$I_{\bar y}=\pi_1(V_{\bar y},\bar \eta)$.

\begin{pr}\label{prKt}
Let the notation be as above.
Let ${\cal F}$ and ${\cal G}$ be
locally constant constructible sheaves
of sets on $V$
of Kummer type
with respect to $D_i,i\in I$
and  put
${\cal H}={\cal H}om(
{\rm pr}_2^*{\cal G},
{\rm pr}_1^*{\cal F})$
on $V\times V$.

1. The base change map
\begin{equation}
\begin{CD}
\tilde \gamma^*\tilde j_*{\cal H}
@>>>
j_{g*}\gamma^*{\cal H}
\end{CD}\label{eqKt}
\end{equation}
is an isomorphism.

2. Let $y\in \widetilde Y_g\setminus V$ be a point
such that $\tilde \delta(y)=\tilde\gamma (y)$.
Let 
$\bar y$ be a geometric point above $y$
and let $I_{\bar y}=\pi_1(V_{\bar y},\bar \eta)$ 
be the inertia group
where $\bar\eta$ is a geometric point
of $V_{\bar y}=V\times_YY_{\bar y}$.
We identify 
$(j_*\delta^*{\cal H})_{\bar y}
=(j_*{\cal H}om({\cal G},{\cal F}))_{\bar y}
=Hom_{I_{\bar y}}({\cal G}_{\bar \eta},
{\cal F}_{\bar \eta})$
and 
$(j_{g*}\gamma^*{\cal H})_{\bar y}
=(j_{g*}{\cal H}om(g^*{\cal G},{\cal F}))_{\bar y}
=Hom_{I_{\bar y}}(g^*{\cal G}_{\bar \eta},
{\cal F}_{\bar \eta})$.
Let $m$ be an integer satisfying the property
in the condition {\rm (Kum)}
for the sheaves ${\cal F}$ and ${\cal G}$.

We define $u_i \in O_{Y,\bar y}^\times$ 
by $g^*(t_i)=u_it_i$ for each $i\in I$.
Then, we have $u_i(\bar y)=1$.
Let $v_i\in O_{Y,\bar y}^\times$
be the unique $m$-th root of $u_i$
such that $v_i(\bar y)=1$ for each $i\in I$.
Let $\bar g\colon \bar \eta\to \bar \eta$
be a map compatible with $g\colon V\to V$
and assume $\bar g^*(\sqrt[m]{t_i})=v_i
\sqrt[m]{t_i}$ for each $i\in I$.
Let $\bar g^*\colon g^*{\cal G}_{\bar \eta}
=\bar g^*({\cal G}_{\bar \eta})
\to
{\cal G}_{\bar \eta}$
be the induced isomorphism.
Then the diagram
\begin{equation}
\begin{CD}
(\tilde j_*{\cal H})_{\tilde\delta(\bar y)}
@>>>
(j_*\delta^*{\cal H})_{\bar y}
=&Hom_{I_{\bar y}}({\cal G}_{\bar \eta},
{\cal F}_{\bar \eta})
\\
@|&@VVV\\
(\tilde j_*{\cal H})_{\tilde\gamma(\bar y)}
@>>>
(j_{g*}\gamma^*{\cal H})_{\bar y}
=&Hom_{I_{\bar y}}(g^*{\cal G}_{\bar \eta},
{\cal F}_{\bar \eta})
\end{CD}
\label{eqKumcan}
\end{equation} 
is commutative,
where the horizontal arrows
are the isomorphisms
{\rm (\ref{eqKt})}
and the right vertical arrow
is induced by the isomorphism
$\bar g^*\colon g^*{\cal G}_{\bar \eta}
\to {\cal G}_{\bar \eta}$.
\end{pr}

{\it Proof.}
1. It is sufficient to show that
the map $(\tilde j_*{\cal H})_{\tilde\gamma(\bar y)}
\to
(j_{g*}\gamma^*{\cal H})_{\bar y}$ is an isomorphism
for each geometric point $\bar y$
of $\widetilde Y_g \setminus V$.
Let $\bar \eta$ be a geometric point
of $V_{\bar y}=V\times_YY_{\bar y}$
where $Y_{\bar y}$ is the strict localization.
Let $I_{\bar y}=\pi_1(V_{\bar y},\bar \eta)$,
$I_{\tilde \gamma(\bar y)}=
\pi_1((V\times V)_{\gamma(\bar y)},
\tilde \gamma(\bar \eta))$
and $I_{g(\bar y)}=
\pi_1(V_{g(\bar y)},g(\bar \eta))$
be the inertia groups.
We regard 
the $I_{\bar y}$-set
${\cal F}_{\bar \eta}$
as an $I_{\tilde \gamma(\bar y)}$-set
by the map
$p_{1*}\colon I_{\tilde \gamma(\bar y)}\to I_{\bar y}$.
We also regard 
the $I_{g(\bar y)}$-set
${\cal G}_{g(\bar \eta)}$
as an $I_{\tilde \gamma(\bar y)}$-set
and an $I_{\bar y}$-set
by the maps
$p_{2*}\colon I_{\tilde \gamma(\bar y)}\to I_{g(\bar y)}$
and by 
$g_*\colon I_{\bar y}\to I_{g(\bar y)}$.
Then, 
the stalks
$(\tilde j_*{\cal H})_{\tilde\gamma(\bar y)}$
and 
$(j_{g*}\gamma^*{\cal H})_{\bar y}$
are naturally identified with
the sets
$Hom_{I_{\tilde \gamma(\bar y)}}
({\cal G}_{g(\bar \eta)},{\cal F}_{\bar \eta})$
and
$Hom_{I_{\bar y}}
({\cal G}_{g(\bar \eta)},{\cal F}_{\bar \eta})$
respectively.
Since $p_{1*}\circ \tilde\gamma_*={\rm id}$
and $p_{2*}\circ \tilde\gamma_*=g_*$, 
it suffices to show that
the images of
$I_{\bar y}$ and of
$I_{\tilde \gamma(\bar y)}$
by the compositions
$$\begin{CD}
I_{\bar y}@>{\tilde \gamma_*}>>
I_{\tilde \gamma(\bar y)}
@>{(p_{1*},p_{2*})}>>
I_{\bar y}\times I_{g(\bar y)}
@>>>
{\rm Aut}({\cal F}_{\bar \eta})
\times
{\rm Aut}({\cal G}_{g(\bar \eta)})
\end{CD}$$
are the same.

Take a basis $t_i$ of $O(-D_i)$
for each $i\in I$ at $\bar y$
and an integer $m$ invertible at $y$
such that
the image of $I_{\bar y}
\to {\rm Aut}({\cal F}_{\bar \eta})$
is a quotient of
the image of the map
$t\colon I_{\bar y}
\to M=\prod_{i\in I}\mu_m$
(\ref{eqiner}).
We also take
a basis $s_i$ of $O(-D_i)$
for each $i$ at $g(y)$
such that
the image of $I_{g(\bar y)}
\to {\rm Aut}({\cal G}_{g(\bar \eta)})$
is a quotient of
the image of
$s\colon I_{g(\bar y)}
\to M=\prod_{i\in I}\mu_m$.
Since $\bar y$ is a geometric point of
$(Y\times Y)^\sim$,
there exists a unit $u_i$
on a neighborhood of $\bar y$
such that $g^*(s_i)=u_it_i$
for each $i\in I$.
Since the $m$-th roots of $u_i$
define an \'etale covering,
we have a commutative diagram
$$\begin{CD}
I_{\tilde \gamma(\bar y)}
@>{(p_{1*},p_{2*})}>>
I_{\bar y}\times I_{g(\bar y)}\\
@VVV @VV{t\times s}V\\
M@>{\Delta}>>M\times M.
\end{CD}$$
Since the map $\tilde \gamma_*\colon 
I_{\bar y}
\to 
I_{\tilde \gamma(\bar y)}$
is a section of
$p_{1*}\colon 
I_{\tilde \gamma(\bar y)}\to I_{\bar y}$,
both
${\rm Image}(I_{\bar y}\to M\times M)$
and
${\rm Image}(I_{\tilde \gamma(\bar y)}\to M\times M)$
are equal to
$\Delta({\rm Image}\ t\colon 
I_{\bar y}\to M)$.
Thus the assertion is proved.

2. The assertion is \'etale local.
We may assume there exists 
a basis $t_i$ of $O(-D_i)$ for each $i\in I$
and ${\cal G}$ is trivialized by 
the \'etale covering
$W=V[T_i,i\in I]/(T_i^m-t_i,i\in I)
\to V$.
Thus, it is reduced to the case
where ${\cal G}$ is the locally constant sheaf
represented by 
the \'etale covering
$W\to V$.

Let $\tilde u_i\in 
\Gamma((Y\times Y)^\sim,O^\times)$
be the unit defined by the equation
$1\otimes t_i=\tilde u_i(t_i\otimes 1)$.
Since $\tilde \gamma(y)=\tilde \delta(y)$,
we have $u_i(y)=\tilde u_i(\tilde \gamma(y))=
\tilde u_i(\tilde \delta(y))=1$.
Hence there exists a unique unit $\tilde v_i\in 
\Gamma((Y\times Y)^\sim_{\bar y},O^\times)$
on the strict localization
such that $\tilde v_i^m=\tilde u_i$ and 
$\tilde v_i(\bar y)=1$.
Then, by 
sending 
$1\otimes T_i$ to $\tilde v_i(T_i\otimes 1)$,
we obtain an isomorphism
$p_2^*{\cal G}\to p_1^*{\cal G}$.
The pull-back of this isomorphism by $\gamma$
defines an isomorphism
$g^*{\cal G}\to {\cal G}$.
The diagram
(\ref{eqKumcan})
is commutative
if we define the right vertical arrow
to be that induced by 
this isomorphism $g^*{\cal G}\to {\cal G}$.

On the other hand,
we have $u_i=\tilde \gamma^*(\tilde u_i)$
and hence $v_i=\tilde \gamma^*(\tilde v_i)$
for each $i\in I$.
By the assumption
$\bar g^*(\sqrt[m]{t_i})=v_i
\sqrt[m]{t_i}$ for each $i\in I$,
the map $\bar g^*\colon g^*{\cal G}_{\bar \eta}
\to {\cal G}_{\bar \eta}$
is also induced by the isomorphism
$g^*{\cal G}\to {\cal G}$
defined above.
Thus the assertion follows.
\qed

We apply Proposition \ref{prKt}
to the following situation.
We consider a cartesian diagram
\begin{equation}
\begin{CD}
V@>>> Y\\
@VfVV @VV{\bar f}V\\
U@>>> X
\end{CD}
\label{eqXYU}
\end{equation}
of normal schemes  
over a perfect field $k$
satisfying the following conditions.
The horizontal arrows are open immersions
with dense images.
The map $f\colon V\to U$ is a finite \'etale Galois
covering of Galois group $G$.
The map $\bar f\colon Y\to X$ is proper
and the action of $G$ on $V$
is extended to an action on $Y$ over $X$.

Let $D_i\subset Y,i\in I$
be a finite family of Cartier divisors of $Y$
such that $V\cap D_i=\emptyset$ for each $i\in I$.
Let ${\cal F}$ be a smooth $E$-sheaf
on $U$ and put ${\cal F}_V=f^*{\cal F}$ on $V$.
We assume that
${\cal F}_V$ is of Kummer type with respect to $D_i,i\in I$.
We put ${\cal H}={\cal H}om(
{\rm pr}_2^*{\cal F}_V,{\rm pr}_1^*{\cal F}_V)$
on $V\times V$.

For $\sigma\in G$,
applying the construction of the diagram (\ref{eqjg})
to the map $\sigma\colon V\to V$,
we obtain an open immersion $j_{\sigma}\colon V\to \widetilde Y_\sigma$
and an immersion $\tilde\gamma_\sigma\colon 
\widetilde Y_\sigma\to (Y\times Y)^\sim$.
We identify
$\delta^*{\cal H}=
{\cal E}nd({\cal F}_V)$
and
$\gamma^*{\cal H}=
{\cal H}om(\sigma^*{\cal F}_V,{\cal F}_V)$.
Let $y\in \widetilde Y_\sigma\setminus V$
be a point satisfying 
$\tilde \gamma_\sigma(y)=\tilde \delta(y)$
and $\bar y$ be a geometric point above $y$.
The base change maps induce isomorphisms
\begin{equation}
\begin{CD}
(j_{\sigma*}{\cal H}om(\sigma^*{\cal F}_V,{\cal F}_V))_{\bar y}
@<<<
\tilde j_*{\cal H}_{\bar y}
@>>>
(j_*{\cal E}nd({\cal F}_V))_{\bar y}
\end{CD}
\label{eqBr}
\end{equation}
by Proposition \ref{prKt}.1.

We recall the definition of the Brauer trace.
Recall that $F$
is the residue field of 
an $\ell$-adic field $E$.
For an automorphism $a$
of an $F$-vector space $M$
of dimension $n$,
the Brauer trace
${\rm Tr}^{Br}(a\colon M)\in E$
is defined as follows.
Let $\det(T-a\colon M)
=(T-\alpha_1)
\cdots(T-\alpha_n)
\in F[T]$
be the eigenpolynomial
of $a$
and let $\tilde \alpha_1,
\ldots,\tilde \alpha_n\in
E^{\rm ur}$
be the Teichm\"uller liftings
of $\alpha_1,
\ldots,\alpha_n$.
Then the Brauer trace
${\rm Tr}^{Br}(a\colon M)$
is defined to be
the sum
$\tilde \alpha_1+
\cdots+\tilde \alpha_n$.

\begin{cor}\label{corKt}
Let the notation be as above.
Assume that the pull-back
${\cal F}_V=f^*{\cal F}$
is of Kummer type with respect to $D_i,i\in I$.
Assume further that the reduction 
$\overline {\cal F}_V$ mod $\lambda$ is constant
and let $\overline M$ denote the corresponding
$F$-representation of $G$.

Let $\sigma\in G$ be an element
of order prime to $\ell$.
We regard 
the canonical map
$\sigma^*\colon \sigma^*{\cal F}_V=
\sigma^*f^*{\cal F}
\to (f\circ \sigma)^*{\cal F}=
f^*{\cal F}={\cal F}_V$
as a section of
$j_{\sigma*}{\cal H}om(\sigma^*{\cal F}_V,{\cal F}_V)$.
Let $y\in \widetilde Y_\sigma\setminus V$ 
be a point such that 
$\tilde \gamma_\sigma (y)=
\tilde \delta(y)$
and $\bar y$ be a geometric point above $y$.

Let $\widetilde{\rm Tr}_{\bar y}
\colon (j_{\sigma*}{\cal H}om(\sigma^*{\cal F}_V,{\cal F}_V))_{\bar y}
\to E$
be the composition of
$$\begin{CD}
(j_{\sigma*}{\cal H}om(\sigma^*{\cal F}_V,{\cal F}_V))_{\bar y}
@>{\rm (\ref{eqBr})}>> 
(j_*{\cal E}nd({\cal F}_V))_{\bar y}
@>{\rm Tr}>> E.
\end{CD}$$
Then, we have
$$\widetilde{\rm Tr}_{\bar y}(\sigma^*)=
{\rm Tr}^{Br}(\sigma\colon \overline M).$$
\end{cor}

{\it Proof.}
We put $\bar x=\bar f(\bar y)$
and take a geometric point $\bar \eta$
of $V_{\bar y}$.
Let $I_{\bar x}=\pi_1(U_{\bar x},\bar \eta)$ be the inertia group.
Then, the element $\sigma\in G$ is in the image of
the natural map $I_{\bar x}\to G$.
By the assumption that 
the order of $\sigma$
is prime to $\ell$,
we may take an inverse image
$\bar \sigma\in I_{\bar x}$
such that the pro-order
of $\bar \sigma$ is prime to $\ell$.
The map $\bar\sigma\colon \bar \eta\to \bar \eta$
is compatible with $\sigma\colon V\to V$.

We take an basis $t_i$ of $O(-D_i)$ for each $i\in I$.
Since $\overline{\cal F}_V$
is constant,
the image of $\pi_1(V,\bar \eta)
\to {\rm Aut}({\cal F}_{V,n,\bar\eta})$
is of order a power of $\ell$
for each $n\ge0$.
Thus, in the condition (Kum) for
${\cal F}_n$,
we may take a power of $\ell$
as an integer $m$.
By the assumption $\tilde\delta(y)=
\tilde\gamma_{\sigma}(y)$,
we have $\sigma(t_i)/t_i=1$ at $y$.
Since $\bar \sigma$ is of pro-order
prime to $\ell$,
for each $i\in I$ and a power $m$ of $\ell$,
we have $\bar\sigma(\sqrt[m]{t_i})/
\sqrt[m]{t_i}=1$ at $\bar y$.
Thus by Proposition \ref{prKt}.2,
the composition
$(j_{\sigma*}{\cal H}om(\sigma^*{\cal F}_V,{\cal F}_V))_{\bar y}
\to (j_*{\cal E}nd({\cal F}_V))_{\bar y}$
is induced by the map $\bar \sigma^*\colon 
\sigma^*{\cal F}_{V,\bar \eta}\to
{\cal F}_{V,\bar \eta}$.
Hence
$\widetilde{\rm Tr}_{\bar y}(\sigma^*)$
is equal to 
${\rm Tr}(\bar \sigma^*\colon {\cal F}_{\bar \eta}).$
Since
the pro-order of $\bar \sigma$
is prime to $\ell$,
the action of $\bar \sigma^*$
on ${\cal F}_{\bar \eta}$
is of finite order prime to $\ell$.
Thus, we have
${\rm Tr}(\bar \sigma^*\colon {\cal F}_{\bar \eta})
=
{\rm Tr}^{Br}(\sigma\colon \overline M).$
\qed

\subsection{Review on the Swan class}

We briefly recall the definition 
of the Swan class \cite{KS}.
First, we recall the definition of the Swan character class.
Assume $k$ is perfect
and we consider a cartesian diagram
\begin{equation}
\begin{CD}
V@>>> Y\\
@VfVV @VV{\bar f}V\\
U@>>> X
\end{CD}
\label{eqXY}
\end{equation}
of schemes  
over $k$.
We assume that the horizontal arrows
are open immersion,
that $\bar f\colon Y\to X$ is proper,
and $f\colon U\to V$ is a finite \'etale
Galois covering of 
smooth schemes of dimension $d$.
Let $G$ be the Galois group
${\rm Gal}(V/U)$.
By alteration \cite{dJ} (cf. \cite{KS} Lemma 3.2.1),
one can construct a diagram
\begin{equation}
\xymatrix{
W\ar[d]_g\ar[rr]&&Z\ar[d]^{\bar g}\ar[ddl]\\
V\ar[dd]_f\ar[rr]&&Y\ar[dd]^{\bar f}\\
&X'\ar[dr]&\\
U\ar[rr]\ar[ru]&&X}
\label{eqXYZ}
\end{equation}
of schemes over $k$
satisfying the following properties.

\begin{enumerate}
\item[(\ref{eqXYZ}.1)]
$U$ is the complement of a Cartier divisor $B$ of $X'$
and the map $X'\to X$ 
is an isomorphism on $U$.
\item[(\ref{eqXYZ}.2)]
$Z$ is smooth over $k$
and $W$ is the complement
of a divisor with simple normal crossings.
\item[(\ref{eqXYZ}.3)]
The map $\bar g\colon Z\to Y$ is proper,
surjective and generically finite.
The squares are cartesian.
\end{enumerate}

We define the log blow-up
$(Z\times Z)'$,
its open subscheme
$(Z\times Z)^\sim=
(W\times W)'$
and the log diagonal map
$Z=\Delta_Z^{\log}\to
(Z\times Z)^\sim=
(W\times W)'$
as in \S2.2.

Let $(X'\times X')'$
be the blow-up of $X'\times X'$
at $B\times B$
and $p_1,p_2\colon (X'\times X')'\to X'$
be the projections.
Let $(X'\times X')^\sim$
be the maximum open subscheme of $(X'\times X')'$
where $p_1^*B=p_2^*B$.
Let $(Z\times_{X'} Z)^\sim$
be the inverse image
of the diagonal 
$X'\to (X'\times X')^\sim$
by the map
$(Z\times Z)^\sim\to 
(X'\times X')^\sim$.

For an element $\sigma\in G$,
let $\Gamma_\sigma$ denote the graph of $\sigma$.
We consider the Gysin map
$(g\times g)^!\colon CH_d(V\times_UV\setminus \Delta_V)
=\bigoplus_{\sigma\neq 1,\sigma\in G}{\mathbb Z}\cdot 
[\Gamma_\sigma]
\to
CH_d(W\times_UW\setminus W\times_VW)$
in the intersection theory \cite{fulton}.
For $\sigma\neq 1$,
the Swan character class
$s_{V/U}(\sigma)\in CH_0(Y\setminus V)_{\mathbb Q}$
is defined as the intersection product
$$s_{V/U}(\sigma)
=-\frac1{[W:V]}\bar g_*(\Gamma,\Delta_Z^{\log})_{
(Z\times Z)^\sim}$$
by taking
a lifting $\Gamma\in 
CH_d((Z\times_{X'}Z)^\sim \setminus W\times_VW)$
of $(g\times g)^!\Gamma_\sigma \in 
CH_d(W\times_UW \setminus W\times_VW)$.
For $\sigma=1$,
$s_{V/U}(\sigma)\in CH_0(Y\setminus V)_{\mathbb Q}$
is defined by requiring
$\sum_{\sigma\in G}s_{V/U}(\sigma)=0$.

The following basic facts are proved in \cite{KS}.

\begin{lm}\label{lmsVU}

1. {\rm (\cite{KS} Proposition 3.2.2)}
The $0$-cycle class $s_{V/U}(\sigma)\in 
CH_0(Y\setminus V)_{\mathbb Q}$
is independent of the choice of
$Z,W$ and $X'$.

Let $h\colon Y'\to Y$
be a proper birational morphism
inducing the identity on $V$.
Then,
we have $h_*s_{V/U}(\sigma)=s_{V/U}(\sigma)$.

2. {\rm (\cite{KS} Lemma 1.1.3, 3.3.2)}
Assume the order of $\sigma$ is not a power of $p$
and let $\overline C\subset (Z\times_{X'}Z)^\sim$
be the closure of
the inverse image 
$C=(g\times g)^{-1}(\Gamma_\sigma)$
of the graph $\Gamma_\sigma
\subset V\times_UV$
of $\sigma$.
Then,
there exists a diagram
{\rm (\ref{eqXYZ})}
such that the intersection
$\overline C\cap \Delta_Z^{\log}$ is empty
and hence
$s_{V/U}(\sigma)=0$.
\end{lm}

We consider an open immersion
$U\to X$
of schemes over $k$.
We assume $U$ is smooth of dimension $d$.
For a smooth $E$-sheaf ${\cal F}$ on $U$,
the naive Swan class is defined as follows.
We take a finite \'etale Galois covering $V\to U$
of Galois group $G$,
trivializing the reduction $\overline{\cal F}$ modulo $\lambda$.
Let $\overline M$
be the $F$-representation of $G$
corresponding to $\overline{\cal F}$.
We take a cartesian diagram (\ref{eqXY})
such that $\bar f\colon Y\to X$ is proper.
Let $G_{(p)}$ denote the subset of
$G$ consisting of elements of order a power of $p$.
Then,
the naive Swan class
${\rm Sw}^{\rm naive}({\cal F})\in 
CH_0(X\setminus U)_E$
is defined by
$${\rm Sw}^{\rm naive}({\cal F})=
\frac1{|G|}\sum_{\sigma\in G_{(p)}}
\bar f_*s_{V/U}(\sigma){\rm Tr}^{Br}(\sigma\colon 
\overline M)$$
where ${\rm Tr}^{Br}$
denotes the Brauer trace.
The naive Swan class
${\rm Sw}^{\rm naive}({\cal F})$
lies in fact in 
$CH_0(X\setminus U)_{E_0}$
where $E_0=E\cap 
{\mathbb Q}(\zeta_{p^\infty})$.
The Swan class
${\rm Sw}({\cal F})
\in 
CH_0(X\setminus U)_{{\mathbb Q}}$
is defined to be
$\frac 1{[E_0:{\mathbb Q}]}
{\rm Tr}_{E_0/{\mathbb Q}}
{\rm Sw}^{\rm naive}({\cal F})$.
They are conjectured to be the same
and known to be equal if
$\dim X\le 2$ \cite{KS} Lemma 4.2.3.2.
It is also proved in the proof 
of Theorem 4.3.14 in loc.\ cit.\
that
$\deg {\rm Sw}({\cal F})
=\deg {\rm Sw}^{\rm naive}({\cal F})$,
if $X$ is proper.

\subsection{The characteristic class and the Swan class}

Now, we state and prove the main result of this section.

\begin{thm}\label{thmSw}
Let $X$ be a scheme  
over a perfect field $k$
and $j\colon U\to X$ be an open immersion.
We assume $U$ is smooth of dimension $d$.
Let ${\cal F}$ be a smooth $E$-sheaf on $U$.

1. 
Let $f\colon V\to U$ be a finite \'etale Galois 
covering trivializing 
the reduction
$\overline{\cal F}$
mod $\lambda$.
Let 
\begin{equation}
\begin{CD}
V@>{j_V}>>Y\\
@VfVV @VV{\bar f}V\\
U@>{j}>>X
\end{CD}
\label{eqXYUV}
\end{equation}
be a cartesian diagram of schemes
over $k$
where $\bar f\colon Y\to X$ is proper.
Let $G={\rm Gal}(V/U)$ be the Galois group 
and $\overline M$ be the 
$F$-representation of
$G$ corresponding to
$\overline{\cal F}$.
Then, if ${\cal F}_V=f^*{\cal F}$ is of Kummer type
with respect to $Y$,
we have
\begin{equation}
C(j_{V!}{\cal F}_V,\overline{\Gamma_\sigma},\sigma^*)
=-s_{V/U}(\sigma){\rm Tr}^{Br}(\sigma\colon \overline M)
\label{eqSwch}
\end{equation}
for $\sigma\in G,\neq 1$.
If the order of $\sigma$ is not
a power of $p$,
the both sides are $0$.

2. Assume ${\cal F}$ is potentially 
of Kummer type with respect to $X$.
Then, we have
\begin{equation}
C(j_!{\cal F})=
{\rm rank}\ {\cal F}\cdot
C(j_!E)-{\rm Sw}^{\rm naive}{\cal F}
\label{eqSw}
\end{equation}
in $H^0(X,{\cal K}_X)$.
\end{thm}

{\it Proof.}
1. 
First, we show that we may assume
the additional conditions 
(\ref{eqXYUV}.1)-(\ref{eqXYUV}.5)
below are satisfied.
By Proposition \ref{prpush} and
Lemma \ref{lmsVU},
we may replace $X$ and $Y$ by
proper modifications.
By replacing $X$
by a blow-up,
we may assume that $U$
is the complement of a
Cartier divisor $B$.
By replacing $Y$
by the closure of the image of
the immersion
$(\sigma)_{\sigma\in G}\colon 
V\to (\prod_{\sigma\in G})_XY$,
we may assume the action of $G$
on $V$ is extended to an action on $Y$.
Let $D_i,i\in I$ be a finite family 
of Cartier divisors of $Y$
with respect to which ${\cal F}_V$
is of Kummer type.
By replacing it by the family
$\sigma(D_i), (i,\sigma)\in I\times G$,
we may assume
the family $D_i,i\in I$
is stable under the $G$-action.
Further, by replacing $Y$
by the blow-up by
the intersections
$D_i\cap D_j$,
we may assume that 
$O(-D_i)+O(-D_j)$
is an invertible ideal
for each $i,j\in I$.
Replacing $X$ and $Y$ by the normalizations,
we may assume $X$ and $Y$ are normal.
Thus we may assume the following conditions
are satisfied.
\begin{enumerate}
\item[(\ref{eqXYUV}.1)]
$U\subset X$ is the complement of a
Cartier divisor $B$.
\item[(\ref{eqXYUV}.2)]
The action of $G$
on $V$ over $U$ is extended to an action on $Y$
over $X$.
\item[(\ref{eqXYUV}.3)]
${\cal F}_V$
is of Kummer type 
with respect to a family $D_i,i\in I$
of Cartier divisors
and the family $D_i,i\in I$
is stable under the $G$-action.
\item[(\ref{eqXYUV}.4)]
$X$ and $Y$ are normal.
\item[(\ref{eqXYUV}.5)]
$O(-D_i)+O(-\sigma^*D_i)$
is an invertible ideal
for each $i\in I$ and $\sigma\in G$.
\end{enumerate}

Let $\overline{(Y\times Y)}'$
be the normalization
of the blow-up of $Y\times Y$
at $D_i\times D_i,i\in I$
as in \S3.1.
Let $(Y\times Y)^\sim
\subset \overline{(Y\times Y)}'$
be the maximum open subscheme
where $p_1^*D_i=p_2^*D_i$ for each $i\in I$
and $\tilde j_V\colon V\times V\to
(Y\times Y)^\sim$ 
be the open immersion.
The closed immersion
$\gamma=(1,\sigma)\colon Y\to Y\times Y$
is uniquely lifted to a closed immersion
$\tilde \gamma\colon Y\to \overline{(Y\times Y)}'$,
by the condition 
(\ref{eqXYUV}.5),
by the universality of blow-up
and the by the assumption that $Y$ is normal.
The intersection 
$\widetilde Y_\sigma=Y\cap (Y\times Y)^\sim$ 
is the maximum open subscheme
where $D_i=\sigma^*D_i$ for each $i\in I$
and $j_\sigma\colon V\to \widetilde Y_\sigma$ 
be the open immersion.

We take an alteration 
\begin{equation}
\begin{CD}
W@>{j_W}>> Z\\
@VgVV@VV{\bar g}V\\
V@>{j_V}>> Y
\end{CD}
\label{eqYZVW}
\end{equation}
satisfying the following properties.
\begin{enumerate}
\item[(\ref{eqYZVW}.1)]
The diagram is cartesian.
\item[(\ref{eqYZVW}.2)]
The map $g\colon Z\to Y$
is proper, surjective and generically finite.
\item[(\ref{eqYZVW}.3)]
$Z$ is smooth over $k$ and $W$
is the complement of a divisor 
$E=\bigcup_{j\in J}E_j$ with
simple normal crossings.
\end{enumerate}

We consider 
the smooth $E$-sheaf
${\cal F}_W=(f\circ g)^*{\cal F}$ on $W$,
the closed subscheme
$C=(g\times g)^{-1}(\Gamma_\sigma)
\subset W\times W$
and the pull-back of the map
$\sigma^*\colon \sigma^*{\cal F}_V\to {\cal F}_V$
on $C$.
We show that they satisfy the assumptions in
Proposition \ref{prtame}.
Let $(Z\times Z)^\sim \subset (Z\times Z)'$
be the log blow-up as in \S2.2.
Since $\overline{\cal F}$
is constant on $V$,
the sheaf ${\cal F}_W$
is tamely ramified along
$Z\setminus W$.
By Lemma \ref{lmtameC},
the subscheme 
$W\times_UW\subset
W\times W$
is non-expanding with respect to $Z$.
Hence its closed subscheme
$C$ is also non-expanding
with respect to $Z$.
Let $\widetilde C$ denote
the closure of $C$ in
$(Z\times Z)^\sim$.
We put ${\cal H}=
{\cal H}om({\rm pr}_2^*{\cal F}_V,
{\rm pr}_1^*{\cal F}_V)$
on $V\times V$.
We consider the closed immersions
$\gamma\colon V\to V\times V$ and 
$\tilde \gamma\colon \widetilde Y_\sigma
\to (Y\times Y)^\sim$.
By Proposition \ref{prKt}
and by the assumption that
${\cal F}_V$ is of Kummer type,
the restriction map 
$\Gamma(\widetilde Y_\sigma,
(\tilde j_{V*}{\cal H})|_{\widetilde Y_\sigma})
\to \Gamma(V,\gamma^*{\cal H})$
is an isomorphism.
By the commutative diagram
$$\begin{CD}
\Gamma(\widetilde Y_\sigma,
(\tilde j_{V*}{\cal H})|_{\widetilde Y_\sigma})
&\quad =\quad & \Gamma(V,\gamma^*{\cal H})&
=Hom_V(\sigma^*{\cal F},{\cal F})
\ni \sigma^*\\
@VVV@VVV&\\
\Gamma(\widetilde C,(\tilde j_{V*}{\cal H})|_{\widetilde C})
&\quad \subset\quad &
\Gamma(C,{\cal H}|_C),&
\end{CD}$$
the image of $\sigma^*$
in $\Gamma(C,{\cal H}|_C)$
lies in the tame part
$\Gamma(\widetilde C,(\tilde j_{V*}{\cal H})|_{\widetilde C})$.
Thus the reductions
${\cal F}_n$ satisfy the assumptions in
Proposition \ref{prtame}.
By applying it to ${\cal F}_n$ 
and taking the limit, we obtain
$$
C(j_{W!}{\cal F}_W,\Gamma_\sigma,\sigma^*)
=
{\rm cl}((g\times g)^!\Gamma_\sigma,
\Delta_W^{\log})_{
(W\times W)'}
\cdot
\widetilde {\rm Tr}\ \delta^{\prime *}\tilde \sigma$$
for $\sigma\in G$.

If the order of $\sigma$
is not a power of $p$,
the intersection
$\widetilde C\cap 
\Delta_Z^{\log}$ is empty
and hence the both sides in
(\ref{eqSwch}) are 0
by Lemma \ref{lmsVU}.2.
If the order of
$\sigma$ is a power of $p$,
we have
$\widetilde {\rm Tr}\ \delta^{\prime *} \tilde \sigma
={\rm Tr}^{Br}(\sigma\colon \overline M)$
by Corollary \ref{corKt}
since ${\cal F}_V$ is assumed of Kummer type.
Thus (\ref{eqSwch}) follows from the definition
of the Swan character class.

2. 
We take a commutative diagram
(\ref{eqXYUV})
as in 1 such that
the pull-back ${\cal F}_V$
is of Kummer type
with respect to $Y$.
Applying Lemma \ref{lmupull}
to the reductions
${\cal F}_n$ and taking the limit, we obtain
$$f^*
(C(j_!{\cal F})-{\rm rank}{\cal F}\cdot
C(j_!E))
=
\sum_{\sigma\in G}(
C(j_{V!}{\cal F}_V,\Gamma_\sigma,\sigma^*)-
{\rm rank}{\cal F}\cdot
C(j_{V!}E,\Gamma_\sigma,1)).$$
By 1,
the right hand side is equal to
$$-\sum_{\sigma\in G_{(p)},\neq 1}
s_{V/U}(\sigma)\cdot 
({\rm Tr}^{Br}(\sigma\colon \overline M)-
{\rm rank}{\cal F})
=-
\sum_{\sigma\in G_{(p)}}
s_{V/U}(\sigma)\cdot 
{\rm Tr}^{Br}(\sigma\colon \overline M).$$
Thus, taking $f_*$, we obtain
\begin{eqnarray*}
|G|(C(j_!{\cal F})-{\rm rank}{\cal F}\cdot
C(j_!E))
&=&
-\sum_{\sigma\in G_{(p)}}
f_*s_{V/U}(\sigma)\cdot 
{\rm Tr}^{Br}(\sigma\colon \overline M)\\
&=&-|G|
{\rm Sw}^{\rm naive}{\cal F}.
\end{eqnarray*}
\qed

We recover the following
main result of \cite{KS}.

\begin{cor}\label{corSw}
{\rm (\cite{KS} Theorem 4.2.9)}
For a smooth $E$-sheaf 
${\cal F}$ on $U$,
we have
\begin{equation}
\chi_c(U,{\cal F})=
{\rm rank}\ {\cal F}
\cdot \chi_c(U,E)-
\deg {\rm Sw}{\cal F}.
\label{eqGOS}
\end{equation}
\end{cor}

{\it Proof.}
Since
$\deg {\rm Sw}{\cal F}=
\deg {\rm Sw}^{\rm naive}{\cal F}$,
it suffices to show
$$\chi_c(U,\overline{\cal F})=
{\rm rank}\ \overline{\cal F}
\cdot \chi_c(U,E)-
\deg {\rm Sw}^{\rm naive}\overline{\cal F}$$
for a smooth $F$-sheaf 
$\overline{\cal F}$ on $U$.
By Brauer induction,
we may assume the rank of 
$\overline{\cal F}$ is 1.
Let $\bar\chi\colon \pi_1(U)^{\rm ab}
\to F^\times$ be the character
corresponding
to $\overline{\cal F}$
and let $\chi\colon \pi_1(U)^{\rm ab}
\to E^\times$ be 
its Teichm\"uller lifting.
Then, the order of $\chi$
is finite
and the corresponding
smooth $E$-sheaf 
${\cal F}_{\chi}$ on $U$
is potentially of Kummer type.
Thus, it suffices to take the trace
of (\ref{eqSw}).
\qed

\section{Characteristic class
of a sheaf of rank 1}

In this section 
we assume $\Lambda$ is a finite
local ${\mathbb Z}_\ell$-algebra.
We fix an inclusion 
${\mu}_{p^\infty}(\Lambda)\to
{\mathbb Z}[\frac1p]/{\mathbb Z}$
and identify the $p$-primary part
$\mu_{p^\infty}(\Lambda)$
as a subgroup of ${\mathbb Z}[\frac1p]/{\mathbb Z}$.
Let $X$ be a smooth scheme over $k$
and $j\colon U\to X$
be the open immersion of the complement
of a divisor $D$ with
simple normal crossings.
For a smooth rank 1 sheaf ${\cal F}$
on $U$,
we will show that the difference
$C(j_!{\cal F})-
{\rm rank}\ {\cal F}\cdot
C(j_!\Lambda)$
is given by the 0-cycle class $c_{\cal F}$
defined by Kato in \cite{Kato2},
assuming ${\cal F}$ is clean with respect to $D$.
The definition of the 0-cycle class $c_{\cal F}$
will be recalled in Definition \ref{dfcF}.2.
As a byproduct,
we obtain a new proof of the
Grothendieck-Ogg-Shafarevich formula
$$\chi_c(U,{\cal F})=
{\rm rank}\ {\cal F}
\cdot \chi_c(U,\Lambda)-
\deg {\rm Sw}{\cal F}.
\leqno{(\ref{eqGOS})}
$$
for curves
without using the Weil formula.

\subsection{Review on ramification
of Artin-Schreier-Witt characters}

We briefly recall the ramification
theory of Artin-Schreier-Witt characters
according to \cite{Kato} \S 3
and \cite{aml} \S 10.
Let $K$ be a complete discrete valuation
field of characteristic $p>0$.

Let $F$ be the residue field of $K$.
Let $\Omega_F=\Omega^1_{F/F^p}$ be 
the $F$-vector space
of differential 1-forms
and put
$$\Omega_F(\log)=
\Omega_F\oplus (F \otimes_{\mathbb Z} K^\times)/
(da -(a\otimes a)\colon a\in O_K,a\neq0).$$
For $a\in K^\times$,
let $d\log a\in \Omega_F(\log)$ 
denote the image of $1\otimes a$.
We have a short exact sequence
\begin{equation}
\begin{CD}
0@>>>
\Omega_F@>>>
\Omega_F(\log)
@>{\rm res}>>
F@>>>0\end{CD}
\label{eqres}
\end{equation}
of $F$-vector space
where the residue map
${\rm res}\colon \Omega_F(\log)\to F$
sends $d\log a$ to ${\rm ord}\ a$.
A choice of prime element
defines a splitting
$F\to \Omega_F(\log)$
of the exact sequence.

A filtration $F_\bullet$ 
on $H^1(K,{\mathbb Q}/{\mathbb Z})$
is defined in \cite{Kato} Definition (2.1)
and is recalled in \cite{aml} \S 10.4.
We have 
$F_{-1}H^1(K,{\mathbb Q}/{\mathbb Z})=0$
and 
$F_0H^1(K,{\mathbb Q}/{\mathbb Z})$
is the unramified part.
For a character $\chi \in 
H^1(K,{\mathbb Q}/{\mathbb Z})$,
the Swan conductor of $\chi$ is defined to 
be the smallest integer $r\ge 0$
such that $\chi \in F^r
H^1(K,{\mathbb Q}/{\mathbb Z})$.
For $r>0$, a canonical injection
\begin{equation}
\begin{CD}
{\rm rsw}_r\colon Gr^F_r
H^1(K,{\mathbb Q}/{\mathbb Z})
@>>> 
\Omega_F(\log)\otimes
m_K^{-r}/m_K^{-r+1}
\end{CD}
\label{eqrsw}
\end{equation}
is defined in Theorem (3.2) (3) \cite{Kato}
and is recalled in Proposition 10.7 \cite{aml}.

\subsection{Blow-up of the diagonal}

Let $X$ be a smooth scheme 
of dimension $d$ over 
a perfect field $k$ and $U\subset X$
be the complement
of a divisor $D=\bigcup_{i\in I}D_i$
with simple normal crossings.
Let $j\colon U\to X$ denote
the open immersion and
let ${\cal F}$ be
a smooth $\Lambda$-sheaf of rank 1
on $U$.

For each irreducible component $D_i$
of $D$,
the sheaf ${\cal F}$
defines a character
of the absolute Galois group
of the local field $K_i$
at the generic point.
By the identification
$\mu_{p^\infty}(\Lambda)
\subset
{\mathbb Z}[\frac1p]/{\mathbb Z}$
fixed at the beginning of this section,
the $p$-primary part 
of the character
$\pi_1(U)^{\rm ab}\to \Lambda^\times$
corresponding to ${\cal F}$
defines a character
$\pi_1(U)^{\rm ab}\to 
{\mathbb Z}[\frac1p]/{\mathbb Z}$.
Hence applying the theory
recalled in the previous subsection,
we define the Swan conductor
${\rm Sw}_i(\cal F)\in {\mathbb N}$
for each $i$.
We define the Swan divisor
$D_{\cal F}$
to be
$\sum_{i\in I}{\rm Sw}_i({\cal F})D_i$.
We decompose
$D=D_t+D_w$
into the sum of
the tame part
$D_t=
\sum_{{\rm Sw}_i({\cal F})=0}D_i$
and the wild part 
$D_w=D_{{\cal F},{\rm red}}
=\sum_{{\rm Sw}_i({\cal F})>0}D_i$.
The refined Swan conductor
${\rm rsw}_i{\cal F}$
defines a non-zero section
of the locally free sheaf
$\Omega^1_{X/k}(\log D)(D_{\cal F})|_{D_w}$
at the generic point
$\xi_i$ of each irreducible component $D_i\subset D_w$.
By \cite{Kato} Theorem (7.1) and Proposition (7.3), 
there exists a unique global section
${\rm rsw}{\cal F}\in
\Gamma(D_w,\Omega^1_{X/k}(\log D)(D_{\cal F})
|_{D_w})$
extending
${\rm rsw}_i{\cal F}$
for each $D_i\subset D_w$.
For a locally free sheaf ${\cal E}$
of finite rank on a scheme $S$,
let $c_*({\cal E})=
\sum_ic_i({\cal E})
\in \bigoplus_i CH^i(S\to S)$
denote the total Chern class
and put $c_*({\cal E})^*=
c_*({\cal E}^*)=
\sum_i(-1)^ic_i({\cal E})$.

\begin{df}\label{dfcF}
Let $X$ be
a smooth scheme of dimension $d$
over a perfect field $k$
and $U\subset X$ be
the complement of a divisor $D$
with simple normal crossings.
Let ${\cal F}$ be
a smooth $\Lambda$-sheaf on $U$
of rank $1$
and ${\rm rsw}{\cal F}
\in \Gamma(D_w,\Omega^1_{X/k}(\log D)(D_{\cal F})
|_{D_w})$
be the refined Swan character
defined on the wild part $D_w$.

1. We say that ${\cal F}$
is clean with respect to $D$
if the section ${\rm rsw}{\cal F}$
is nowhere vanishing on $D_w$.

2. If ${\cal F}$
is clean with respect to $D$,
we define a $0$-cycle class
$c_{\cal F}\in CH_0(D_w)$
by
$$c_{\cal F}=
\{c_*(\Omega^1_{X/k}(\log D))^*\cap
(1+D_{\cal F})^{-1}
\cap [D_{\cal F}]\}_{\dim 0}.$$
\end{df}

The refined Swan character
may be also regarded as an $O_{D_w}$-linear map
${\rm rsw}{\cal F}\colon O_{D_w}(-D_{\cal F})\to
\Omega^1_{X/k}(\log D)|_{D_w}$.
Then, 
we have
$$c_{\cal F}=
(-1)^{d-1}c_{d-1}({\rm Coker}
({\rm rsw}{\cal F}\colon 
O_{D_w}(-D_{\cal F})\to
\Omega^1_{X/k}(\log D)|_{D_w}))\cap 
[D_{\cal F}].$$
In \cite{KS},
the equality
$c_{\cal F}={\rm Sw}({\cal F})$
is proved if $d\le 2$
in Theorem 5.1.5
and is conjectured in general
in Conjecture 5.1.1.

Let $f\colon (X\times X)'\to X\times X$
be the log blow-up 
with respect to $D\subset X$
and let $X\to (X\times X)'$
be the log diagonal as in \S2.2.
Recall that
$(U\times X)',(X\times U)'\subset (X\times X)'$ are
the complement
of the proper transforms
of $D\times X$ and of $X\times D$
respectively.
In this section, we
write $(U\times U)'$
for the log product
$(X\times X)^\sim=
(U\times X)'\cap
(X\times U)'$.

We regard $D_{\cal F}\subset X$ 
as a closed subscheme of 
$(X\times X)'$ by
the log diagonal map
$X\to (X\times X)'$
and let
$\bar g\colon (X\times X)''
\to (X\times X)'$
be the blow-up at
$D_{\cal F}$.
Let $(U\times U)'',
(U\times X)'',(X\times U)''
\subset 
(X\times X)''$
be the inverse image of
$(U\times U)',
(U\times X)',(X\times U)'$
respectively.
Let $\Delta_i\subset (X\times X)'$
be the exceptional divisor
above $D_i\times D_i$
for each $i\in I$.
Let 
$(U\times U)'''
\subset 
(U\times U)''$
be the complement of the union of
the proper transforms of
$\Delta_i$ for $D_i\subset D_w$.
We consider a commutative diagram
$$\xymatrix{
& (X\times U)''\ar[dl]_{j''_2}\ar[dd]& 
&(U\times U)''
\ar[ll]_{k''_1}\ar[dl]_{k''_2}\ar[dd]_{g}&
(U\times U)'''\ar[l]_{\tilde j''}\\
(X\times X)''\ar[dd]_{\bar g}& &
(U\times X)''\ar[dd]\ar[ll]_{\!\!\!\!\!\! j''_1}&&&\\
& (X\times U)' \ar[dl]_{j'_2}& 
&(U\times U)'\ar[ll]_{\!\!\! \!\!\! k'_1}\ar[dl]_{k'_2}&&\\
(X\times X)'\ar[dd]_f& &
(U\times X)'\ar[ll]_{\!\!\!\!\!\! j'_1}&&&\\
& X\times U \ar[dl]_{j_2}\ar[uu]&
& U\times U.\ar[ll]_{\!\!\! \!\!\! k_1}\ar[dl]_{k_2}
\ar[uu]_{\tilde j'}\ar[uuuur]_{\tilde j'''}&\\
X\times X& &U\times X\ar[ll]_{j_1}\ar[uu]&&&}$$
The vertical down arrows are blowing-ups
and the other arrows are open immersions.
The  faces of the upper cube
and the four faces consisting of
open immersions of the lower cube
are cartesian.

The projections
$p_1,p_2\colon (U\times U)'''\to X$
are smooth of dimension $d$
and hence $(U\times U)'''$
is smooth of dimension $2d$
over $k$.
The log diagonal $X\to (U\times U)'$
is uniquely lifted to
a closed immersion $X\to (U\times U)'''$.
We identify the exceptional divisor
of the blow-up
$g\colon (U\times U)''\to
(U\times U)'$
with the ${\mathbf P}^d$-bundle
${\mathbf P}(\Omega^1_{X/k}(\log D)|_{D_w}
\oplus N_{D_{\cal F}/X})
={\mathbf P}(
\Omega^1_{X/k}(\log D)(D_{\cal F})|_{D_w}
\oplus O_{D_w})$
over $D_w$.
Under this identification,
the inverse image
$(U\times U)'''\times_{(U\times U)}D_w
=(U\times U)'''\setminus
(U\times U)$
of $D_w$ is identified with the 
${\mathbf A}^d$-bundle
${\mathbf V}(\Omega^1_{X/k}(\log D)(D_{\cal F})|_{D_w})$
over $D_w$.
Here 
${\mathbf P}({\cal E})=
{\bf Proj}(S^\bullet{\cal E})$
and
${\mathbf A}({\cal E})=
{\bf Spec}(S^\bullet{\cal E})$
for a locally free sheaf ${\cal E}$.

\begin{pr}\label{prrk1}
Let $X$ be
a smooth scheme of dimension $d$ over 
a perfect field $k$
and $U\subset X$ be
the complement of a divisor
$D$ with simple normal crossings.
Let ${\cal F}$ be
a smooth $\Lambda$-sheaf on $U$
of rank $1$
and $D_{\cal F}=
{\rm Sw}({\cal F})$
be the Swan divisor
as above.
We put ${\cal H}_0={\cal H}om
({\rm pr}_2^*{\cal F},
{\rm pr}_1^*{\cal F})$
on $U\times U$.
Then, we have the following.

1. The $\Lambda$-sheaf
${\cal H}_0'''=
\tilde j'''_*{\cal H}_0$
is a smooth $\Lambda$-sheaf of rank 1
on $(U\times U)'''$.

2. The restriction of
the smooth $\Lambda$-sheaf
${\cal H}_0'''$
on the complement
$(U\times U)'''\setminus
(U\times U)$
is the Artin-Schreier
sheaf defined by
the minus of the refined Swan character
${\rm rsw}{\cal F}\in 
\Gamma(D_w,
\Omega^1_{X/k}(\log D)(D_{\cal F})
|_{D_w})$
regarded as a linear form on
$(U\times U)'''\setminus
(U\times U)=
{\mathbf V}(\Omega^1_{X/k}(\log D)(D_{\cal F})|_{D_w})$.
\end{pr}

First, we verify it at
the generic point
of each irreducible component
of $(U\times U)'''\setminus
(U\times U)=
{\mathbf V}(\Omega^1_{X/k}(\log D)(D_{\cal F})|_{D_w})$.
Let $D_i$ be an irreducible component of $D_w$
and $\xi_i$ be the generic point of $D_i$.
Let $K_i$ be the local field of $X$ at $\xi_i$ and
let $\chi_i\colon G_{K_i}^{\rm ab}\to \Lambda^\times$
be the character defined by the sheaf $\cal F$.
The function field $F_i$ of $D_i$
is the residue field of $K_i$.
Put $n_i={\rm Sw}_{K_i}\chi_i>0$
and let ${\rm rsw} \chi_i \in
\Omega_{F_i}(\log)\otimes 
m_{K_i}^{-n_i}/m_{K_i}^{-n_i+1}$
denote the refined Swan character of 
the $p$-primary part of $\chi_i$.
Let $\eta_i$ be the generic point of the divisor
$(U\times U)'''\times_
{(U\times U)'}D_i$,
let $L_i$ be the local field
at $\eta_i$
and $E_i$ be the residue field.
Let
$\varphi_i\colon G_{L_i}^{\rm ab}\to \Lambda^\times$
be the character defined by the sheaf ${\cal H}_0$.
We have an inclusion
$\Omega_{F_i}(\log)\otimes 
m_{K_i}^{-n_i}/m_{K_i}^{-n_i+1}=
\Omega^1_{X/k}(\log D)(D_{\cal F})_{\xi_i}
\subset E_i$.

The following lemma is a consequence of
\cite{aml} Proposition 13.6.

\begin{lm}\label{lmrk1}
Let the notation be as above
and $n_i={\rm Sw}_{K_i}\chi_i>0$.
Then we have the following. 

1. The character 
$\varphi_i\colon G_{L_i}^{\rm ab}\to \Lambda^\times$
is unramified
and is of order $p$.

2. The character $\varphi_i$ regarded as an element of
$H^1(L_i,{\mathbb Z}/p{\mathbb Z})$
is the image of
the minus of the refined Swan character 
${\rm rsw}\chi_i
\in \Omega_{F_i}(\log)\otimes 
m^{-n_i}_{\xi_i}/
m^{-n_i+1}_{\xi_i}$
by the canonical map
$$\begin{CD}
\Omega_{F_i}(\log)\otimes 
m^{-n_i}_{\xi_i}/
m^{-n_i+1}_{\xi_i}
\subset
E_i@>>>
H^1(E_i,{\mathbb Z}/p{\mathbb Z})
\subset 
H^1(L_i,{\mathbb Z}/p{\mathbb Z})
\end{CD}$$
of Artin-Schreier theory.
\end{lm}

\begin{cor}\label{corlmrk1}
The smooth sheaf ${\cal H}_0$ on $U\times U$
is ramified along the component
$\Delta_i$ of the  
exceptional divisor
$(U\times U)'\setminus 
(U\times U)$ over $D_i\times D_i$
if $n_i={\rm Sw}_{K_i}\chi_i>0$.
\end{cor}

{\it Proof.}
We prove it by contradiction.
Assume the sheaf ${\cal H}_0$
is unramified along $\Delta_i$.
Then the sheaf ${\cal H}_0$
is extended to a smooth sheaf 
on a neighborhood of
the generic point $\xi_i$ of $D_i
=\Delta_i\cap \Delta_X^{\log}$.
Then, the character
$\varphi_i
\in H^1(L_i,{\mathbb Z}/p{\mathbb Z})$
is defined by the pull-back
of a character
in $G_{\kappa(\xi_i)}^{\rm ab}=
G_{F_i}^{\rm ab}\to \Lambda^\times$.
This contradicts to
Lemma \ref{lmrk1}.2.
\qed

\medskip

{\it Proof of Proposition \ref{prrk1}.}
1. It follows from the assertion 1 in
Lemma \ref{lmrk1}
and the purity of branch locus. 

2. It follows from the assertion 2
in Lemma \ref{lmrk1}.
\qed

\begin{cor}\label{cork1}
Let the notation be as in Proposition \ref{prrk1}.
We put ${\cal H}_0'=
\tilde j'_*{\cal H}_0$.
Then, we have the following.

1. The canonical map
$j_{1!}Rk_{2*}{\cal H}_0
=Rj_{2*}k_{1!}{\cal H}_0
\to
Rf_*j'_{1!}Rk'_{2*}{\cal H}'_0
=Rf_*Rj'_{2*}k'_{1!}{\cal H}'_0$
is an isomorphism.

2. If ${\cal F}$
is clean with respect to $D$,
the identity of ${\cal H}_0$
is extended to an isomorphism
${\cal H}_0'\to
Rg_*\tilde j''_!
{\cal H}_0'''.$
\end{cor}

{\it Proof.}
1. It suffices to show that
the restrictions of
$Rf_*j'_{1!}Rk'_{2*}{\cal H}'_0$
and $Rf_*Rj'_{2*}k'_{1!}{\cal H}'_0$
are 0 on $D\times D$.
Since the assertion is \'etale local on $X$,
we may assume the following conditions
are satisfied:
We have $X=X_1\times X_2$,
$U=U_1\times U_2$
and ${\cal F}={\cal F}_1\boxtimes {\cal F}_2$.
The open subschemes
$U_1\subset X_1$ and
$U_2\subset X_2$ are
the complement of divisors
$D_1\subset X_1$ and
$D_2\subset X_2$
with simple normal crossings.
The rank 1 sheaf
${\cal F}_1$ on $U_1$ is
tamely ramified along $D_1$
and 
${\cal F}_2$ on $U_2$ is
wildly ramified along
each component of $D_2$.

Then, by the K\"unneth formula,
it suffices to show
the cases where $X=X_1$ and $X=X_2$
respectively.
The case ${\cal F}$ is tamely ramified
is proved in Lemma \ref{lmFP}.
Assume 
${\cal F}$ is
wildly ramified along
each component of $D$.
Then, by Corollary \ref{corlmrk1},
the sheaf ${\cal H}_0$
is ramified along each component
of the exceptional divisor
$(U\times U)'\setminus (U\times U)$.
Thus we have ${\cal H}_0'=
\tilde j'_!{\cal H}_0$
and the assertion follows.

2. Let $g_0\colon T=(U\times U)'''\setminus (U\times U)
\to D_w$ denote the restriction
of $g\colon (U\times U)''\to (U\times U)'$.
It suffices to show
$Rg_{0!}({\cal H}_0'''|_T)=0$.
By Proposition \ref{prrk1}.2,
the restriction ${\cal H}_0'''|_T$
is the Artin-Schreier sheaf
defined by the refined Swan conductor
${\rm rsw}{\cal F}$
on the ${\mathbf A}^d$-bundle
$T={\mathbf V}(\Omega^1_{X/k}(\log D)(D_{\cal F})|_{D_w})$.
By the assumption that 
${\rm rsw}{\cal F}$
is nowhere vanishing,
the restriction of
${\cal H}_0'''$ on
every geometric fiber
is an Artin-Schreier sheaf
on ${\mathbf A}^d$ defined by
a non-trivial linear form.
Thus the assertion follows.
\qed

Now, we are ready to prove the
main result of this section.

\begin{thm}\label{thmrk1}
Let $X$ be
a smooth scheme of dimension $d$
over a perfect field $k$
and $U\subset X$ be
the complement of a divisor
with simple normal crossings.
Let ${\cal F}$ be
a smooth $\Lambda$-sheaf on $U$
of rank $1$.
If ${\cal F}$ is clean with respect to
$D$, we have
$$C(j_!{\cal F})=
C(j_!\Lambda)-c_{\cal F}$$
in $H^{2d}(X, \Lambda(d))$.
\end{thm}

{\it Proof.}
We keep the notation in Proposition \ref{prrk1}.
We put ${\cal H}_0=
{\cal H}om({\rm pr}_2^*{\cal F},
{\rm pr}_1^*{\cal F})$
and ${\cal H}_0^*=
{\cal H}om({\rm pr}_1^*{\cal F},
{\rm pr}_2^*{\cal F})$
on $U\times U$.
By Proposition loc.\ cit.,
the sheaves ${\cal H}_0'''=
\tilde j'''_*{\cal H}_0$
and ${\cal H}_0^{\prime\prime\prime*}=
\tilde j'''_*{\cal H}_0^*$
on $(U\times U)'''$
are smooth of rank 1.
The natural pairing
${\cal H}_0\otimes {\cal H}_0^*
\to \Lambda$ on $U\times U$
is extended uniquely
to a pairing
${\cal H}_0'''
\otimes
{\cal H}_0^{\prime\prime\prime*}
\to \Lambda$
on $(U\times U)'''$.

Further we put
$\overline{\cal H}
=j_{1!}Rk_{2*}
{\cal H}_0(d)[2d]$
and 
$\overline{\cal H}^*
=j_{2!}Rk_{1*}
{\cal H}_0^*(d)[2d]$
on $X\times X$
and 
$\overline{\cal H}''
=j''_{1!}Rk''_{2*}
\tilde j''_!
{\cal H}_0'''(d)[2d]$
and 
$\overline{\cal H}^{\prime\prime*}
=j''_{2!}Rk''_{1*}
\tilde j''_!
{\cal H}_0^{\prime\prime\prime*}(d)[2d]$
on $(X\times X)''$.
We define 
$\tilde 1\in H^0_X((X\times X)'',
\overline{\cal H}'')=
\Gamma(X,j_*{\cal E}nd({\cal F}))
=\Gamma(U,{\cal E}nd({\cal F}))$
to be the identity.
We also define
$\tilde 1^*\in H^0_X((X\times X)'',
\overline{\cal H}^{\prime\prime *})$
to be the identity.
By Corollary \ref{cork1},
we obtain canonical isomorphisms
$\overline{\cal H}
\to R(f\circ \bar g)_*
\overline{\cal H}''$
and 
$\overline{\cal H}^*
\to R(f\circ \bar g)_*
\overline{\cal H}^{\prime\prime*}$.
The pairing
${\cal H}_0'''
\otimes
{\cal H}_0^{\prime\prime\prime*}
\to \Lambda$
on $(U\times U)'''$
induces a pairing
$\overline{\cal H}''
\otimes
\overline{\cal H}^{\prime\prime*}
=(j'_1\circ k'_2\circ \tilde j'')_!({\cal H}_0'''
\otimes
{\cal H}_0^{\prime\prime\prime*}(2d)[4d])
\to {\cal K}_{(X\times X)''}$
on $(X\times X)''$.

Thus the assumptions 
in Proposition \ref{prpink} are satisfied.
Since both $(f\circ \bar g)_*\tilde 1$
and $(f\circ \bar g)_*\tilde 1^*$
are equal to $j_!1=1$,
we obtain 
$C(j_!{\cal F})=
\langle \tilde 1,
\tilde 1^*\rangle$.
By the compatibility
(\ref{eqcupint}),
the right hand side is further equal to
${\rm cl}(X,X)_{(X\times X)''}$.
Since
the conormal sheaf
$N_{X/(X\times X)''}$
is isomorphic to
$\Omega_{X/k}(\log D)(D_{\cal F})$,
the self-intersection product
$(X,X)_{(X\times X)''}$
is equal to the 0-cycle class
$$(-1)^dc_d(\Omega_{X/k}(\log D)(D_{\cal F}))\cap [X]
=\{c_*(\Omega_{X/k}(\log D))(1-D_{\cal F})^{-1}\cap [X]\}^*_{\dim 0}.$$
Thus the assertion follows.
\qed

As in Section 3,
let $E$ denote a finite
extension of ${\mathbb Q}_\ell$.

\begin{cor}\label{corrk1}
Let $X$ be a scheme over 
a perfect field $k$
and $U\subset X$ be a dense open subscheme
smooth over $k$.
Let ${\cal F}$
be a smooth $E$-sheaf
of rank $1$ on $U$
and let $j\colon U\to X$ be the open immersion.
If $\dim U\le 2$, we have
$$C(j_!{\cal F})=
C(j_!E)-{\rm Sw}{\cal F}$$
in $H^0(X, {\cal K}_X)$.
\end{cor}

{\it Proof.}
By Proposition \ref{prpush} and
Lemma \ref{lmsVU},
we may replace $X$ by
a proper modification.
Hence, by \cite{Kato2} Theorem 4.1,
we may assume that $X$ is smooth,
that $U\subset X$ is
the complement of a divisor
with simple normal crossings
and that
${\cal F}$ is clean with respect to the boundary.
Then the assertion follows from
Theorem \ref{thmrk1}
applied to the reductions ${\cal F}_n$
and the equality
${\rm Sw}({\cal F})=c_{\cal F}$
(\cite{KS} Theorem 5.1.5).
\qed

\noindent{\bf Remark.}
By the argument in 
the proof of Corollaries \ref{corSw},
\ref{corrk1},
we obtain
the Grothendieck-Ogg-Shafarevich formula
$$\chi_c(U_{\bar k},{\cal F})=
{\rm rank}{\cal F}
\cdot \chi_c(U_{\bar k}, {\mathbb Q}_\ell)-
\deg {\rm Sw}{\cal F}
\leqno{\rm(\ref{eqGOS})}$$
for a smooth $E$-sheaf 
${\cal F}$ if $\dim U\le 2$.
This proof does not use the Weil formula.

\section{Localized characteristic class.}

Let $X$ be a smooth scheme over
$k$ and $S\subset X$ be a closed subscheme
such that the complement $U$
is dense.
Let ${\cal F}$ be an object
of $D_{\rm ctf}(X)$ and assume that
the cohomology sheaves ${\cal H}^q{\cal F}$
are smooth on $U$.
In this section,
we will define a localization
$C^0_S({\cal F})
\in H^0_S(X,{\cal K}_X)$
of the characteristic class
as a cohomology class
supported on $S$.

\subsection{Discrepancy of
the map (\ref{FoxG}).}
Let 
$i\colon Z\to X$
be a closed immersion of schemes
over $k$.
For an object ${\cal F}$
of $D^-(Y,\Lambda)$,
we have a canonical map
$i^*{\cal F}\otimes^L Ri^!\Lambda
\to Ri^!{\cal F}$
(\ref{FoxG}).
We construct a functor
$$\Delta_i\colon D^-(X,\Lambda)
\to D^-(Z,\Lambda)$$
fitting in a distinguished triangle
\begin{equation}
\to i^*{\cal F}\otimes^LRi^!\Lambda
\to Ri^!{\cal F}
\to \Delta_i {\cal F}\to.
\label{eqphitri}
\end{equation}

Let $P\to X$
be a conservative family of
geometric points.
For a $\Lambda$-sheaf
${\cal F}$ on $X$,
let ${\cal F}\to {\cal C}({\cal F})$
be the Godement
resolution defined by $P$
(\cite{4XVII} 4.2.2).
A canonical map
${\cal F}\otimes i_*i^!{\cal C}(\Lambda)
\to i_*i^!{\cal C}({\cal F})$
is defined as follows.
Let $U\to X$ be an etale morphism
and $s\in {\cal F}(U)$
be a local section.
Then, by regarding $s$
as a map $\Lambda_U\to
{\cal F}|_U$,
we obtain a map
$s\colon i_*i^!{\cal C}(\Lambda)|_U
\to i_*i^!{\cal C}({\cal F})|_U$
by functoriality.
In other words,
we have a map
${\cal F}(U)\otimes
i_*i^!{\cal C}(\Lambda)|_U
\to i_*i^!{\cal C}({\cal F})|_U$.
They define a map
${\cal F}\otimes
i_*i^!{\cal C}(\Lambda)
\to i_*i^!{\cal C}({\cal F})$
or equivalently
$i^*{\cal F}\otimes
i^!{\cal C}(\Lambda)
\to i^!{\cal C}({\cal F})$.

Let ${\cal K}$
be an object of
$D^-(X,\Lambda)$.
Let ${\cal K}_1\to {\cal K}$
be the Cartan-Eilenberg flat resolution.
The map defined above
induces a map
$a\colon i^*{\cal K}_1\otimes
i^!{\cal C}(\Lambda)
\to i^!{\cal C}({\cal K}_1)$
of double complexes.
Let $N>0$ be an integer
greater than the cohomological dimension
of $i^!$.
Let $\tau''_{<N}$
be the partial canonical truncation
with respect to the 
Godement degree.
Then we define 
$$\Delta_i({\cal K})=
{\rm Cone}\left(
\textstyle{\int}(i^*{\cal K}_1\otimes
i^!\tau''_{<N}{\cal C}(\Lambda))
\to 
\textstyle{\int}i^!\tau''_{<N}{\cal C}({\cal K}_1)
\right)$$
to be the mapping cone
of the map of simple complexes
associated to the partial truncation
of the map of double 
complexes.

\begin{lm}\label{lmphi}
For ${\cal K}\in
D^-(X,\Lambda)$,
we have a distinguished
triangle
$$
\to i^*{\cal F}\otimes^LRi^!\Lambda
\to Ri^!{\cal F}
\to \Delta_i {\cal F}\to.
\leqno{\rm(\ref{eqphitri})}$$
\end{lm}
{\it Proof.}
The simple complexes
associated to
$i^*{\cal K}_1\otimes
i^!\tau''_{<N}{\cal C}(\Lambda)$
and
$i^!\tau''_{<N}{\cal C}({\cal K}_1)$
compute
$i^*{\cal F}\otimes^LRi^!\Lambda$
and
$Ri^!{\cal F}$
respectively.
Further the map $a$
defines
the canonical map
$i^*{\cal F}\otimes^LRi^!\Lambda
\to Ri^!{\cal F}$.
Hence, the assertion follows.
\qed

\begin{cor}\label{corphi}
If the cohomology sheaves
${\cal H}^q{\cal K}$
are smooth
on a open subscheme
$U\subset X$,
the complex $\Delta_i{\cal K}$
is acyclic on $Z\cap U$.
\end{cor}
{\it Proof.}
The canonical map
$i^*{\cal F}\otimes^LRi^!\Lambda
\to Ri^!{\cal F}$
is an isomorphism on $U\cap Z$.
Hence, the assertion follows
from the distinguished triangle
(\ref{eqphitri}).
\qed

\subsection{Localized characteristic class}
Let $X$ be a scheme 
over $k$
and $\delta\colon X\to X\times X$
be the diagonal map.
Let ${\cal F}$ be
an object of $D_{\rm ctf}(X)$
and put
${\cal H}=
R{\cal H}om({\rm pr}_1^*{\cal F},
{\rm pr}_2^!{\cal F})$
on $X\times X$.
Recall that
the characteristic class
$C({\cal F})
\in H^0(X,{\cal K}_X)$
is defined by the composition
$\delta_*\Lambda
\to {\cal H}
\to \delta_*{\cal K}_X$
by Proposition \ref{prKS}.
We consider
\begin{equation}
\Delta_{\delta}\delta_*\Lambda
\to \Delta_{\delta}{\cal H}
\to \Delta_{\delta}\delta_*{\cal K}_X.
\label{eqphi}
\end{equation}

Assume the cohomology sheaves
${\cal H}^q{\cal F}$ are
smooth on the complement $U$
of a closed subscheme $S\subset X$.
Then, the cohomology sheaves of
${\cal H}$
are smooth on $U\times U$
and hence
the complex $\Delta_{\delta}{\cal H}$ on 
the diagonal $X$
is acyclic on $U$
by Corollary \ref{corphi}.
Hence the composition
of the canonical map
$\Lambda_X=
R\delta^!\delta_*\Lambda_X\to 
\Delta_{\delta}\delta_*\Lambda$
with the composition of 
(\ref{eqphi})
defines a cohomology class
$C_S^0({\cal F})\in
H^0_S(X, 
\Delta_{\delta}\delta_*{\cal K}_X)$.

Further assume $X$ is smooth 
of dimension $d$.
Since the diagonal
$\delta\colon X\to X\times X$
is a section of
the smooth projection
${\rm pr}_2\colon X\times X\to X$,
we have an isomorphism
${\cal K}_X\otimes R\delta^!\Lambda
\to \Lambda_X$.
Hence, by Corollary \ref{corphi}, 
we obtain a distinguished triangle
\begin{equation}
\to
\Lambda
\to 
{\cal K}_X
\to
\Delta_\delta\delta_*{\cal K}_X
\to
\label{eqpure}
\end{equation}
and hence
an exact sequence
$$
H^0_S(X, \Lambda)
\to 
H^0_S(X, {\cal K}_X)
\to
H^0_S(X, \Delta_\delta\delta_*{\cal K}_X)
\to 
H^1_S(X, \Lambda)
.$$
We also assume that
$U=X\setminus S$ is
dense in $X$.
Then,
we have
$H^0_S(X, \Lambda)
=H^1_S(X, \Lambda)=0$
and hence
the canonical map
$H^0_S(X, {\cal K}_X)
\to
H^0_S(X, \Delta_\delta\delta_*{\cal K}_X)$
is an isomorphism.
Thus, in this case,
the class
$C_S^0({\cal F})\in
H^0_S(X, \Delta_\delta\delta_*{\cal K}_X)$
gives an element of
$H^0_S(X, {\cal K}_X)$.

\begin{df}
Let $X$ be a smooth scheme over
$k$ and $S\subset X$ be a closed subscheme
such that the complement $U$
is dense.
Let ${\cal F}$ be an object
of $D_{\rm ctf}(X)$ and assume that
the cohomology sheaves ${\cal H}^q{\cal F}$
are smooth on $U$.
We call 
$C_S^0({\cal F})\in
H^0_S(X, {\cal K}_X)$
the localized characteristic class
of ${\cal F}$.
\end{df}

\begin{lm}\label{lmladd}
Let $X$ be a smooth scheme over $k$
and $({\cal F},F_\bullet)$
be an object of 
$DF_{\rm ctf}(X)$.
Assume that
the cohomology sheaves
${\cal H}^qGr^F_p{\cal F}$ are
smooth on the complement
of a closed subscheme $S\subset X$
for each $p\in {\mathbb Z}$.
Then, we have
$$C^0_S({\cal F})=\sum_qC^0_S(Gr^F_q{\cal F})$$
in $H^0_S(X,{\cal K}_X)$.
\end{lm}

{\it Proof.}
The proof is similar
to that of Lemma \ref{lmadd}.
The map $\delta_*\Lambda
\to {\cal H}=
R{\cal H}om(
{\rm pr}_2^*{\cal F},
{\rm pr}_1^!{\cal F})$
is induced by a map
$\delta_*\Lambda
\to F_0{\cal H}=
F_0R{\cal H}om(
{\rm pr}_2^*{\cal F},
{\rm pr}_1^!{\cal F})$.
We have a canonical isomorphism
$Gr^F_0{\cal H}
\to
\bigoplus_q R{\cal H}om(
{\rm pr}_2^*Gr^F_q{\cal F},
{\rm pr}_1^!Gr^F_q{\cal F})$.
The map
${\cal H}\to \delta^!{\cal K}_X$
is induced by
the sum of
the maps
${\cal H}_q=R{\cal H}om(
{\rm pr}_2^*Gr^F_q{\cal F},
{\rm pr}_1^!Gr^F_q{\cal F})
\to \delta^!{\cal K}_X$.
Thus the composition
$\delta_*\Lambda
\to F_0{\cal H}
\to {\cal H}\to
\delta^!{\cal K}_X$
is equal to the sum
of the compositions
$\delta_*\Lambda
\to F_0{\cal H}
\to {\cal H}_q\to
\delta^!{\cal K}_X$.
Since the cohomology sheaves
of $F_0{\cal H}$
are smooth on $U\times U$,
the assertion follows.
\qed

We show that
the localized characteristic class
refines the characteristic class of ${\cal F}$.

\begin{pr}\label{prloc}
Let $X$ be a smooth scheme over
$k$ and $S\subset X$ be a closed subscheme
such that the complement $U$
is dense.
Let ${\cal F}$ be an object
of $D_{\rm ctf}(X)$ and assume that
the cohomology sheaves ${\cal H}^q{\cal F}$
are smooth on $U$.
Then the image of
the localized
characteristic class
$C_S^0({\cal F})\in
H^0_S(X, {\cal K}_X)$
by the canonical map
$H^0_S(X, {\cal K}_X)\to
H^0(X, {\cal K}_X)$ is
equal to the difference
$C({\cal F})
-{\rm rank}{\cal F}|_U\cdot
C(\Lambda)$.
\end{pr}

{\it Proof.}
Let $j\colon U\to X$
be the open immersion.
By applying Lemma \ref{lmladd}
to the filtration defined
by the subcomplex
$j_!{\cal F}_U\subset
{\cal F}$,
it is reduced to the cases
where
${\cal F}_U=0$
and 
${\cal F}=j_!{\cal F}_U$
respectively.

First, we show the case
${\cal F}_U=0$.
In this case,
the complex ${\cal H}$
is supported on $S\times S$
and the localized class
$C_S^0({\cal F})\in 
H^0_S(X,{\cal K}_X)$
is defined by
the composition
$\delta_*\Lambda
\to {\cal H}\to 
\delta_*{\cal K}_X
\to \delta_*\Delta_\delta
\delta_*{\cal K}_X$.
Thus the assertion 
follows in this case.

We show the case
${\cal F}=j_!{\cal F}_U$.
First, we consider
the case ${\cal F}=j_!\Lambda$.
By Corollary \ref{corphi},
we have $C_S^0(\Lambda)=0$.
Hence we have
$C_S^0(j_!\Lambda)=-
C_S^0(i_*\Lambda)$
by Lemma \ref{lmladd}
and its image is
$-C(i_*\Lambda)=
C(j_!\Lambda)-
C(\Lambda)$.
Thus the assertion is proved in this case.

We prove the general case
${\cal F}=j_!{\cal F}_U$.
Since it has been proved for
$j_!\Lambda_U$,
it suffices to show that
the image of
$C_S^0(j_!{\cal F})
-{\rm rank}\ {\cal F}\cdot C_S^0(j_!\Lambda)$
is equal to
$C(j_!{\cal F})
-{\rm rank}\ {\cal F}\cdot C(j_!\Lambda)$.
To show this, we define a variant
$C_{S!}^0(j_!{\cal F}_U)\in
H^0_S(X, \Delta_\delta\delta_*j_!{\cal K}_U)$
of the difference
$C_!(j_!{\cal F}_U)
-{\rm rank}\ {\cal F}
\cdot C_!(j_!{\cal F}_U)
\in
H^0(X, j_!{\cal K}_U)$.
Recall that
the refinement $C_!(j_!{\cal F}_U)$
of $C(j_!{\cal F}_U)$
is defined in
Definition \ref{dfccc}
by the composition of
$\delta_*\Lambda
\to {\cal H}
\to \delta_*j_!{\cal K}_U$.
Applying $\Delta_\delta$
to this,
we obtain a cohomology class
$C_{S!}^0(j_!{\cal F}_U)\in
H^0_S(X,\Delta_\delta\delta_*j_!{\cal K}_U)$.
Since its image in
$H^0_S(X,\Delta_\delta\delta_*{\cal K}_X)
\overset\sim\gets
H^0_S(X,{\cal K}_X)$
is 
$C_S^0(j_!{\cal F}_U)$,
it suffices to
show the following.

\begin{lm}\label{lmloc}
1. The difference
$C_{S!}^0(j_!{\cal F}_U)
-{\rm rank}{\cal F}\cdot
C_{S!}^0(j_!\Lambda)$
is in the image of
the injection
$H^0_S(X, j_!{\cal K}_U)
\to H^0_S(X,\Delta_\delta\delta_*j_!{\cal K}_U)$.

2. The image
of its inverse image 
in $H^0_S(X, j_!{\cal K}_U)$
by the map
$H^0_S(X, j_!{\cal K}_U)
\to 
H^0(X,j_!{\cal K}_U)$ is
equal to the difference
$C_!(j_!{\cal F})
-{\rm rank}{\cal F}\cdot
C_!(j_!\Lambda)$.
\end{lm}
{\it Proof.}
We may assume $X$ and hence
$U$ is connected.
We may also assume $S$ is non-empty.
We consider a commutative diagram
$$\begin{CD}
@. @.@.
H^0(U,\Lambda_U)\\
@.@. @. @VVV\\
H^0_S(X,j_!\Lambda_U)@>>>
H^0_S(X,j_!{\cal K}_U)@>a>>
H^0_S(X,\Delta_\delta\delta_*j_!{\cal K}_U)
@>{\partial}>>
H^1_S(X,j_!\Lambda_U)\\
@.@VVV @VVV @VVV\\
H^0(X,j_!\Lambda_U)@>>>
H^0(X,j_!{\cal K}_U)@>b>>
H^0(X,\Delta_\delta\delta_*j_!{\cal K}_U)@>>>
H^1(X,j_!\Lambda_U)
\end{CD}$$
of exact sequences.
The horizontal arrows
are induced by the 
distinguished triangle 
$\to 
j_!\Lambda_U\to
j_!{\cal K}_U
\to\Delta_\delta\delta_*j_!{\cal K}_U
\to$, similar to
(\ref{eqpure}).
The classes
$C_{S!}^0(j_!{\cal F})$
and
$C_{S!}^0(j_!\Lambda)$
lie in third term
$H^0_S(X,\Delta_\delta\delta_*j_!{\cal K}_U)$
on the upper line
and
$C_!(j_!{\cal F})$
and
$C_!(j_!\Lambda)$
lie in the second term
$H^0(X,j_!{\cal K}_U)$
on the lower line.
By the definition,
they have the same images
in the third term
$H^0(X,\Delta_\delta\delta_*j_!{\cal K}_U)$
on the lower line.

Since
$H^0_S(X,j_!\Lambda_U)=
H^0(X,j_!\Lambda_U)=0$,
the horizontal arrows
$a$ and $b$
are injective.
Thus, 
it suffices to prove
\begin{equation}
\partial C_{S!}^0(j_!{\cal F})=
\text{ rank }{\cal F}\cdot
\partial C_{S!}^0(j_!\Lambda)
\label{eqbound}
\end{equation}
in
$H^1_S(X,j_!\Lambda_U)$.

By the exact sequence
$0=H^0_S(X,\Lambda_X)
\to 
H^0_S(X,\Lambda_S)
\to 
H^1_S(X,j_!\Lambda_U)
\to 
H^1_S(X,\Lambda_X)$ $=0$,
the map
$H^0(S,\Lambda_S)
=
H^0_S(X,\Lambda_S)
\to 
H^1_S(X,j_!\Lambda_U)$
is an isomorphism.
Thus, we may regard the equality
(\ref{eqbound})
as an equality in $H^0(S,\Lambda)$.

By the commutative diagram,
the boundary
$\partial C_{S!}^0(j_!{\cal F})$
is in the image of the injection
$H^0(X,\Lambda_X)=
H^0(U,\Lambda_U)\to
H^1_S(X,j_!\Lambda_U)=
H^0(S,\Lambda_S).$
Thus, 
we may further regard the equality
(\ref{eqbound})
as an equality in $H^0(X,\Lambda_X)$.

We show the equality
(\ref{eqbound})
in $H^0(X,\Lambda_X)$,
by reducing it to
the case where ${\cal F}$
is unramified along $S$.
Take a closed point $x\in U$
and put $S'=S\amalg\{x\}$
and $U'=X\setminus S'$.
Let $j'\colon U'\to X$
be the open immersion
and put ${\cal F}'={\cal F}|_{U'}$.
Then, in the direct sum decomposition
$H^0(S',\Lambda)=
H^0(S,\Lambda)
\oplus
H^0(x,\Lambda)$,
we have
$C^0_{S'!}(j'_!{\cal F}')=
(\partial C_{S!}^0(j_!{\cal F}),
\partial C^0_{x!}(j'_!{\cal F}))
\in
H^0(S,\Lambda)
\oplus
H^0(x,\Lambda)$.
It lies in the diagonal image
of $H^0(X,\Lambda)$.
Since the construction
is \'etale local,
we have
$\partial C^0_{x!}(j'_!{\cal F})
=\text{ rank }{\cal F}\cdot
\partial C^0_{x!}(j'_!\Lambda)$.
Thus we also have
$\partial C_{S!}^0(j_!{\cal F})
=\text{ rank }{\cal F}\cdot
\partial C_{S!}^0(j_!\Lambda)$
and the assertion is proved.
\qed

\begin{tabular}{l@{\qquad}l}
Ahmed Abbes&                     Takeshi Saito\\
LAGA,&                           Department of Mathematical 
Sciences,\\
Universit\'e de Paris 13,&               University of Tokyo,\\
99, Av. J.-B. Cl\'ement,& Tokyo 153-8914 Japan\\
93430 Villetaneuse France&            \\
abbes@math.univ-paris13.fr&       t-saito@ms.u-tokyo.ac.jp
\end{tabular}

\end{document}